\documentclass[11pt]{amsart}
\usepackage{amssymb}\usepackage{amsmath, amsthm, amsopn, amsfonts}
\usepackage[dvips]{graphicx}\usepackage{graphpap}
\newtheorem{theoreme}{Theorem}

\newtheorem{proposition}{Proposition}
\newtheorem{lemme}[proposition]{Lemma}
\newtheorem{corollaire}[proposition]{Corollary}

\newtheorem{remarque}[proposition]{Remark}

\numberwithin{equation}{section}
\numberwithin{proposition}{section}

 \def\11{{\rm 1~\hspace{-1.4ex}l} }
\def\R{\mathbb R}\def\C{\mathbb C}\def\Z{\mathbb Z}\def\N{\mathbb N}\def\E{\mathbb E}

\begin{document}
\title[Invariant measures for Benjamin-Ono equation]
{ Construction of a Gibbs measure associated to the periodic Benjamin-Ono equation }
\author{N. Tzvetkov}
\address{D\'epartement de Math\'ematiques, Universit\'e Lille I, 59 655
Villeneuve d'Ascq Cedex, France}
\email{nikolay.tzvetkov@math.univ-lille1.fr}
\begin{abstract} 
We define a finite Borel measure of Gibbs type, supported by the Sobolev spaces of negative
indexes on the circle. The measure can be seen as a limit of finite dimensional
measures. These finite dimensional measures are invariant by the ODE's which correspond
to the projection of the Benjamin-Ono equation, posed on the circle, on the
first $N$, $N\geq 1$ modes in the trigonometric bases.
\end{abstract}
\maketitle
%
\section{Introduction, preliminaries and statement of the main result}
Let us denote by $S^1$ the circle identified with $\R/ (2\pi\Z)$. For $u\in
{\mathcal D}'(S^1)$ a distribution on $S^1$, we define its Fourier coefficients as
$\hat{u}(n)\equiv (2\pi)^{-1}u(\exp(-inx))$, $n\in\Z$. Then, we have the Fourier expansion of $u$ (cf. \cite{S}),
$$u=\sum_{n\in\Z}\hat{u}(n)e^{inx}\quad {\rm in}\quad {\mathcal D}'(S^1)\,.$$
We say that $u$ is real valued, if $u=\bar{u}$, where $\bar{u}$ is defined as
$$\bar{u}(\varphi)=\overline{u(\bar{\varphi})},\quad \forall \varphi\in
C^{\infty}(S^1;\C). $$ 
We also have that $u$ is real valued iff its Fourier coefficients satisfy, 
$$\hat{u}(n)=\overline{\hat{u}(-n)},\quad \forall n\in\Z\,.$$
For $s\in\R$, we denote by $H^{s}(S^1)$ the complex Hilbert space of
distributions on $S^1$ equipped with the scalar product
$\langle\cdot,\cdot\rangle_{s}$, defined by
\begin{equation}\label{sp}
\langle u,v\rangle_{s}=2\pi\,\sum_{n\in\Z}\langle
n\rangle^{2s}\hat{u}(n)\overline{\hat{v}(n)}\,,
\end{equation}
where $\langle n\rangle\equiv (1+n^2)^{1/2}$. 
For $s=0$, $H^s(S^1)=L^2(S^1)$ and for $s\geq 0$, the space  $H^{s}(S^1)$ contains integrable functions
on the circle, while for $s<0$ the elements of $H^{s}(S^1)$ are not induced by
integrable functions via the canonical identification. Denote by
$H^{s}_{0}(S^1)$ the subset of $H^{s}(S^1)$ defined as
$$H^{s}_{0}(S^1)\equiv \{ u\in H^{s}(S^1)\,:\, \hat{u}(0)=0,\quad
\hat{u}(n)=\overline{\hat{u}(-n)},\quad \forall\,n\in\Z^{\star} \}.$$
Notice that the elements of $H^{s}_{0}(S^1)$ are real valued distributions.
We have that $H^{s}_{0}(S^1)$ endowed with the scalar product (\ref{sp}) is a
real Hilbert space.
For $s\geq 0$, the space $H^{s}_{0}(S^1)$ contains the real valued functions of $H^{s}(S^1)$ with zero mean value.
\\
 
Consider the Cauchy problem for the Benjamin-Ono equation, posed on $S^1$,
\begin{equation}\label{BO}
(\partial_{t}+H\partial^{2}_{x})u+\partial_{x}(u^{2})=0,\quad u|_{t=0}=u_{0}\in H^s_{0}(S^1) \end{equation}
for some $s\in\R$.
In (\ref{BO}), $H:H^s_{0}(S^1)\longrightarrow H^s_{0}(S^1)$ denotes the
Hilbert transform defined for 
$w\in H^{s}_{0}(S^1)$ as
$$
w\longmapsto -i\sum_{n\in\Z^{\star}}{\rm sign}(n)\hat{w}(n)e^{inx}\,,
$$
i.e
$$\widehat{Hw}(n)\equiv -i\,{\rm sign}(n)\hat{w}(n),\quad n\in\Z^{\star},
\quad\widehat{Hw}(0)\equiv 0\,. $$
Considering solutions of (\ref{BO}) in the space $H^{s}_{0}(S^1)$ seems
reasonable since by a {\it formal integration} of the equation the mean value
of $u$ is preserved. If $s<0$ the expression $u^2$ is a priori not defined and
the interpretation of the nonlinear term in (\ref{BO}) requires to be done carefully.
For $s\geq 0$, it follows from the work of Molinet \cite{Molinet} that
(\ref{BO}) has a well-defined global in time dynamics in the phase space $H^s_{0}(S^1)$.
\\

Recall that the Benjamin-Ono equation is an asymptotic model derived from the
Euler equation for the propagation of internal long waves (see \cite{Benjamin}).
\\

The goal of this paper is to construct a weighted Wiener measure, of Gibbs
type associated to (\ref{BO}). This construction is in the spirit of the work by
Lebowitz-Rose-Speer \cite{LRS} for the nonlinear Schr\"odinger equation.
As we will see, in the context of (\ref{BO}) the construction requires more
involved probabilistic arguments compared to  \cite{LRS}.
\\

We fix for the remaining part of this paper a positive number $\sigma$. 
The Gibbs type measure we construct will be a finite Borel measure on
$H^{-\sigma}_{0}(S^1)$. For an integer $N\geq 1$, we consider the finite
dimensional sub-space of $H^{-\sigma}_{0}(S^1)$ defined as follows
$$
E_{N}\equiv\, \Big(u\in H^{-\sigma}_{0}(S^1)\,:\, \hat{u}(n)=0,\,\, |n|>N\Big).
$$
Notice that the elements of $E_N$ are real valued $C^{\infty}(S^1)$ functions
and we may identify $E_N$ with $\R^{2N}$ by specifying a bases of $E_N$.
A canonical bases of $E_N$ is formed by $\cos(nx)$, $\sin(nx)$, $1\leq n\leq N$.
One can also equip $E_{N}$ with a canonical measure induced by the mapping
from $\R^{2N}$ to $E_{N}$ defined as follows
\begin{equation}\label{map}
(a_1,\cdots,a_N,b_1,\cdots,b_N)\longmapsto \sum_{n=1}^{N}
\Big(a_{n}\cos(nx)+b_{n}\sin(nx)\Big).
\end{equation}
Let us denote by $S_{N}$ the Dirichlet projector defined for $u\in {\mathcal D}'(S^1)$ as
$$
S_{N}(u)\equiv \sum_{|n|\leq N}\hat{u}(n)e^{inx}\,.
$$
Notice that if $u\in H^{-\sigma}_{0}(S^1)$ then $S_{N}(u)\in E_{N}$. Let us
consider the following ordinary differential equation with phase space $E_{N}$
\begin{equation}\label{BO_N}
(\partial_{t}+H\partial^{2}_{x})u_{N}+S_{N}\big(\partial_{x}(u_{N}^{2})\big)=0,\quad u_{N}|_{t=0}=u_{0}\in E_{N}. 
\end{equation}
Let us decompose $u_{N}(t,x)$ in the canonical bases as
$$
u_{N}(t,x)=\sum_{n=1}^{N}
\Big(a_{n}(t)\cos(nx)+b_{n}(t)\sin(nx)\Big),\quad a_{n}(t),b_{n}(t)\in \R.
$$
Then, if we set
$$
c_{n}(t)\equiv \frac{1}{2}(a_{n}(t)-ib_{n}(t))
$$
we can write
$$
u_{N}(t,x)=\sum_{0<|n|\leq N}c_{n}(t)e^{inx},\quad c_{n}(t)=\overline{c_{-n}(t)}.
$$
Thus (\ref{BO_N}) is an ODE for the coefficients $c_{n}(t)$, $0<|n|\leq N$.
More precisely for $0<|n|\leq N$,
\begin{equation}\label{BO_Nbis}
\dot{c}_{n}(t)=-i\,{\rm sign}(n) n^{2}c_{n}(t)-in
\sum_{\stackrel{ 0<|n_1|\leq N,0<|n_2|\leq N}{ n=n_1+n_2}}
c_{n_1}(t)c_{n_2}(t)\,,\quad c_{n}(0)=\widehat{u_0}(n).
\end{equation}
Observe that the equation for $n$ is the complex conjugate of the equation for
$-n$ and thus (\ref{BO_Nbis}) is a system of $N$ ordinary differential
equation for 
$$c(t)\equiv(c_1(t),\cdots,c_{N}(t))\in \C^{N}$$ which can be written in the form
$\dot{c}=P(c)$ with $P$ a polynomial of $c,\bar{c}$ of degree $2$
(equivalently one may write an ODE of similar type for $(a_n,b_n)$). 
Thus we can apply the Cauchy-Lipschitz theorem for ODE's to  (\ref{BO_Nbis}) and deduce that for
every real valued $u_{0}\in E_{N}$ there exists a unique local in time
solutions of (\ref{BO_Nbis}) on a small time interval. Moreover, either the
solution is global in time or there exists $T\neq 0$ such that 
\begin{equation}\label{blowup}
\lim_{t\rightarrow T}\max_{0<|n|\leq N}|c_{n}(t)|=\infty.
\end{equation}
Since integrations by parts give
$$
\int_{S^1}(H\partial^{2}_{x}(u_{N}))u_{N}=-\int_{S^1}(H\partial_{x}(u_{N}))(\partial_{x}(u_{N}))=0
$$
and
$$\int_{S^1}S_{N}\big(\partial_{x}(u_{N}^{2})\big)u_{N}=\int_{S^1}\partial_{x}(u_{N}^{2})u_{N}
=\frac{2}{3}\int_{S^1}\partial_{x}(u_{N}^{3})=0,
$$
by multiplying (\ref{BO_N}) by $u_{N}$, we obtain that
$$
\partial_{t}\big(\int_{S^1}u_{N}^2(t,x)dx\big)=0.
$$
Thus the local solutions of  (\ref{BO_N}) satisfy
\begin{equation}\label{L2}
\sum_{0<|n|\leq N}|c_{n}(t)|^{2}=\frac{1}{2\pi}\|u_{N}(t,\cdot)\|^{2}_{L^2(S^1)}=\frac{1}{2\pi}\|u_0\|^{2}_{L^2(S^1)}\,.
\end{equation}
Therefore (\ref{blowup}) is excluded and 
thus we obtain that for every $u_{0}\in E_{N}$ the ODE (\ref{BO_N}) has a unique global in
time solution.
\\

The problem (\ref{BO_N}) is a Hamiltonian ODE resulting from the Hamiltonian
$F$ defined by
$$
F(u_{N})\equiv-\frac{1}{2}\int_{S^1}\big(|D_{x}|^{\frac{1}{2}}(u_{N})\big)^{2}-\frac{1}{3}\int_{S^1}u_{N}^{3},
$$
where for $u\in H^{s}(S^1)$ the operator $|D_{x}|^{\frac{1}{2}}$ is
defined as Fourier multiplier by
$$
\widehat{|D_{x}|^{\frac{1}{2}}u}(n)\equiv |n|^{\frac{1}{2}}\hat{u}(n),\quad n\in\Z.
$$
Notice that $H^s_{0}(S^1)$ is invariant under the action of $|D_{x}|^{\frac{1}{2}}$.
In addition, we have that 
$$
|D_{x}|^{\frac{1}{2}}\,\circ\,|D_{x}|^{\frac{1}{2}}=H\partial_x
$$
and for real valued $u,v\in C^{\infty}(S^1)$,
$$
\int_{S^1}\big(|D_{x}|^{\frac{1}{2}}u\big)(x)\,v(x)dx=\int_{S^1}u(x)\, \big(|D_{x}|^{\frac{1}{2}}v\big)(x)dx\,.
$$
We can write (\ref{BO_N}) as
$$
\partial_{t}u_{N}=\frac{d}{dx}\nabla F(u_{N})
$$
where $\nabla$ is the $L^2$ gradient on $E_{N}$. Therefore the Hamiltonian $F$
is also conserved by the flow of (\ref{BO_N}).
Let us give a direct proof of this fact.
We can write (\ref{BO_N}) as
\begin{equation}\label{proj}
\partial_{t}u_{N}+\partial_{x}(H\partial_{x}u_{N}+S_{N}(u_{N}^{2}))=0.
\end{equation}
Multiplying the last equation by $H\partial_{x}u_{N}+S_{N}(u_{N}^{2})$ and
integrating over $S^1$, we get
$$
\int_{S^1}(\partial_{t}u_{N})\big(H\partial_{x}u_{N}+S_{N}(u_{N}^{2})\big)=0
$$
and thus
$$
\frac{1}{2}\partial_{t}\Big(\int_{S^1}\big(|D_{x}|^{\frac{1}{2}}(u_{N})\big)^{2}\Big)+
\int_{S^1}(\partial_{t}u_{N})S_{N}(u_{N}^2)=0\,.
$$
On the other hand, using that $\partial_{t}u_{N}\in E_{N}$, we get
$$
\int_{S^1}(\partial_{t}u_{N})S_{N}(u_{N}^2)=
\int_{S^1}(\partial_{t}u_{N})(u_{N}^2)=
\frac{1}{3}\partial_{t}(\int_{S^1}u_{N}^{3}).
$$
Therefore $\partial_{t}(F(u_{N}(t,\cdot))=0$ which implies the Hamiltonian
conservation for the solutions of (\ref{BO_N}).
\\

Let us now observe that  (\ref{BO_Nbis}) can be written in the coordinates
$a^{N}=(a_1,\cdots,a_{N})$, $b^{N}=(b_1,\cdots,b_{N})$ as
\begin{equation}\label{ab}
\partial_{t}a^{N}=-J_{N}\frac{\partial{F}}{\partial b^N},\quad
\partial_{t}b^{N}=J_{N}\frac{\partial{F}}{\partial a^N},
\end{equation}
where
$J_{N}=-\frac{1}{\pi}{\rm diag}(1,2,\cdots,N)$ and
$$
F=F(a^{N},b^{N})=
F\Big(
\sum_{n=1}^{N}
\big(a_{n}\cos(nx)+b_{n}\sin(nx)\big)
\Big).
$$
Indeed, the projection of (\ref{proj}) on the mode $\cos(nx)$ is
$$
\pi\dot{a}_{n}+\pi n^2 b_n+\int_{S^1}\partial_{x}(S_{N}(u_{N}^{2}))\cos(nx)dx=0
$$
and we may write
\begin{eqnarray*}
\int_{S^1}\partial_{x}(S_{N}(u_{N}(x)^{2}))\cos(nx)d & = &
n\int_{S^1}S_{N}(u_N(x)^2)\sin(nx)dx
\\
& = &
n\int_{S^1}u_N(x)^2\sin(nx)dx
\\ 
& = & 
n\frac{\partial}{\partial b_n}\Big(
\frac{1}{3}\int_{S^1}u_{N}^3
\Big).
\end{eqnarray*}
On the other hand
$$
\int_{S^1}\big(|D_{x}|^{\frac{1}{2}}u_N\big)^{2}(x)dx
=
\pi\sum_{n=1}^{N}n(a_n^2+b_n^2)
$$
and thus
$$
\pi n^2 b_n=n
\frac{\partial}{\partial b_n}\Big(
\frac{1}{2}
\int_{S^1}
\big(|D_{x}|^{\frac{1}{2}}u_N(x)\big)^{2}dx
\Big).
$$
Therefore the projection of (\ref{proj}) on the mode $\cos(nx)$ can be written
as
$$
\pi\dot{a}_{n}=\frac{n}{\pi} \frac{\partial}{\partial b_n}\Big(
F\Big(
\sum_{n=1}^{N}
\big(a_{n}\cos(nx)+b_{n}\sin(nx)\big)
\Big)
\Big).
$$
Similarly
$$
\pi\dot{b}_{n}=-\frac{n}{\pi} \frac{\partial}{\partial a_n}\Big(
F\Big(
\sum_{n=1}^{N}
\big(a_{n}\cos(nx)+b_{n}\sin(nx)\big)
\Big)
\Big)
$$
and thus, the equation (\ref{proj}) may indeed be written in the form (\ref{ab}).
Since
\begin{equation}\label{div}
\frac{\partial}{\partial a^{N}}\Big(-J_{N}\frac{\partial{F}}{\partial b^N}\Big)
+
\frac{\partial}{\partial b^{N}}\Big(J_{N}\frac{\partial{F}}{\partial a^N}\Big)=0
\end{equation}
the Liouville theorem for divergence free vector fields (cf. e.g. \cite{Zh}) applies to 
(\ref{ab}), and thus to (\ref{BO_N}) too.
More precisely, if we denote by $\Phi_{N}(t):E_{N}\rightarrow E_{N}$, $t\in\R$ the flow of
(\ref{BO_N}) then it follows from the Liouville theorem that the Lebesgue
measure $\lambda_{N}$ on $E_{N}$ is invariant by the flow of (\ref{BO_N}). 
Namely, for every measurable set $A\subset E_{N}$ and
every $t\in\R$ one has $\lambda_{N}(A)=\lambda_{N}(\Phi_{N}(A))$.
Since $F$ is a conserved quantity for  (\ref{BO_N}), we also have that for
every $\beta\in\R$ the Gibbs measure $\exp(\beta F(u_{N}))d\lambda_{N}(u_{N})$ 
is also invariant by the flow of (\ref{BO_N}).
Moreover since the $L^2$ norm of $u_{N}$ is also a conserved quantity, we
have that for every real constant $c_{N}$ and every measurable function
$\chi_{N}:\R\rightarrow\R$ the measure 
\begin{equation}\label{measure}
c_{N}\chi_{N}(\|u_{N}\|_{L^2(S^1)})\exp(\beta F(u_{N}))d\lambda_{N}(u_{N})
\end{equation}
is also conserved by the flow of (\ref{BO_N}).
We are going to show that for a suitable choice of $c_{N}$ and $\chi_{N}$ the
measures (\ref{measure}), extended to $H^{-\sigma}_{0}(S^1)$, tend to a
limit measure which is a finite Borel measure on $H^{-\sigma}_{0}(S^1)$,
absolutely continuous with respect to a Wiener measure on  $H^{-\sigma}_{0}(S^1)$
induced by a Gaussian process.
\\

Recall that we identify the Lebesgue measure on $E_N$ as the image measure
under the map (\ref{map}) from $\R^{2N}$ to $E_N$. Let us next consider the
measure $d\tilde{\theta}_{N}$ defined as
$$
d\tilde{\theta}_{N}\equiv
e^{-\pi\sum_{n=1}^{N}n(a_n^{2}+b_{n}^{2})}\prod_{n=1}^{N}da_n db_n\,.
$$
Notice that $\tilde{\theta}_{N}(\R^{2N})=(N!)^{-1}$. We then consider
the probability measure 
$$
d\theta_{N}\equiv N!\,d\tilde{\theta}_{N}.
$$
We still denote by $\theta_{N}$ the measure on $E_N$ induced from $\theta_{N}$ by the mapping
(\ref{map}). 
\\

Let us fix a family $h_{n},l_{n}\in {\mathcal N}(0,1)$, $n=1,2,\cdots$
of independent identically distributed  standard real valued Gaussian
variables on a probability space $(\Omega,{\mathcal A},p)$.
Let us observe that the measure $\theta_{N}$ is the
distribution of the $E_N$ valued random variable defined as
$$
\varphi_{N}(\omega,x)=
\sum_{n=1}^{N}
\Big(
\tilde{h}_{n}(\omega)\cos(nx)+\tilde{l}_{n}(\omega)\sin(nx)
\Big),
$$
where $\tilde{h}_{n}, \tilde{l}_{n}\in {\mathcal N}(0,1/\sqrt{2\pi n})$ are
independent identically distributed real Gaussian random variables on
$(\Omega,{\mathcal A},p)$. Thus we may assume that
$
\tilde{h}_{n}(\omega)=(2\pi n)^{-\frac{1}{2}}h_{n}(\omega)
$
and
$
\tilde{l}_{n}(\omega)=(2\pi n)^{-\frac{1}{2}}l_{n}(\omega),
$
where 
$h_{n},l_{n}\in {\mathcal N}(0,1)$
are the fixed standard real valued Gaussians.
Therefore, if we set
$$
g_{n}(\omega)\equiv\frac{1}{\sqrt{2}}(h_{n}(\omega)-il_{n}(\omega))
$$
then 
$
(g_{n}(\omega))_{n=1}^{N}
$
is a sequence of standard independent identically distributed  complex Gaussians and
\begin{equation*}
\varphi_{N}(\omega,x)=\sum_{0<|n|\leq N}\frac{g_{n}(\omega)}{2\sqrt{\pi |n|}}\,
e^{inx},\quad g_{n}(\omega)=\overline{g_{-n}(\omega)}.
\end{equation*}
Let us denote by $L^2(\Omega;H^{-\sigma}_{0}(S^1))$ the Banach space of
$H^{-\sigma}_{0}(S^1))$ valued functions on $\Omega$ (the integration of such
functions being understood in the sense of Bochner integrals).
Clearly $(\varphi_{N})$ is a Cauchy sequence in
$L^2(\Omega;H^{-\sigma}_{0}(S^1))$ and hence we can define
\begin{equation}\label{stat}
\varphi(\omega,x)=\sum_{n\neq 0}\frac{g_{n}(\omega)}{2\sqrt{\pi |n|}}\,
e^{inx},\quad g_{n}(\omega)=\overline{g_{-n}(\omega)}.
\end{equation}
as an element of $L^2(\Omega;H^{-\sigma}_{0}(S^1))$. In particular
$\varphi(\omega,\cdot)\in H^{-\sigma}_{0}(S^1)$ almost surely
and the map $\omega\mapsto \varphi(\omega,\cdot)$ is measurable from $(\Omega,{\mathcal A})$
to $(H^{-\sigma}_{0}(S^1),{\mathcal B})$, where ${\mathcal B}$ denotes the
Borel sigma algebra of $H^{-\sigma}_{0}(S^1)$.
Thus $\varphi(\omega,x)$ defines a measure $\theta$ on
$(H^{-\sigma}_{0}(S^1),{\mathcal B})$ as follows : if $A\in {\mathcal B}$
then 
$
\theta(A)\equiv p(\omega\,:\, \varphi(\omega,\cdot)\in A).
$
Let $\chi_{R}:\R\rightarrow[0,1]$ be a continuous function with compact
support such that $\chi_{R}(x)=1$ for $|x|\leq R$.
Define the measure $d\mu_{N}$ on $E_N$ as
$$
d\mu_{N}(u_N)=\chi_{R}\Big(\|u_{N}\|_{L^2(S^1)}^{2}-\alpha_{N}\Big)e^{-\frac{2}{3}\int_{S^1}u_{N}(x)^{3}dx}
d\theta_{N}(u_{N}),
$$
where
$$
\alpha_{N}\equiv
\sum_{n=1}^{N}\frac{1}{n}=\E\Big(\|\varphi_{N}(\omega,\cdot)\|^{2}_{L^2(S^1)}\Big)\,.
$$
Notice that $\alpha_{N}$ diverges as $\log(N)$ for $N\gg 1$. 
Observe that in the coordinates $a_{n},b_{n}$ given by (\ref{map}) the measure
$\mu_{N}$ reads
$$
N!\,\,\chi_{R}\Big(\|u_{N}\|_{L^2(S^1)}^{2}-\alpha_{N}\Big)e^{2F(u_{N})}\prod_{n=1}^{N}da_n
db_n,
$$
with
$$
u_{N}=\sum_{n=1}^{N}\big(a_{n}\cos(nx)+b_{n}\sin(nx)\big).
$$
From the above discussion (see (\ref{div})) the measure $\prod_{n=1}^{N}da_n
db_n$ is invariant and since $F$ and the $L^2$ norm are conserved under the
flow of (\ref{BO_N}), we obtain that $d\mu_{N}$ is invariant under the flow of (\ref{BO_N}).
\\

Observe that if $A\in {\mathcal B}$ is a Borel set of $H^{-\sigma}_{0}(S^1)$
then $A\cap E_{N}$ is a Borel set of $E_{N}$
(indeed, this is clear for cylindrical sets $A$ and then can be
extends to all $A\in {\mathcal B}$ using that ${\mathcal B}$ is the minimal
sigma algebra containing all cylindrical sets).
We then define the measure
$\rho_{N}$ which is the natural extension of $\mu_{N}$ to
$(H^{-\sigma}_{0}(S^1),{\mathcal B})$. More precisely for every $A\in {\mathcal B}$ which 
is a Borel set of $H^{\sigma}_{0}(S^1)$, we set
$$
\rho_{N}(A)\equiv \mu_{N}(A\cap E_{N})\,.
$$
We now can state the main result of this paper.
\begin{theoreme}\label{thm}
The sequence
\begin{equation}\label{cinm}
\chi_{R}\Big(\|S_{N}(u)\|_{L^2(S^1)}^{2}-\alpha_{N}\Big)e^{-\frac{2}{3}\int_{S^1}(S_{N}u)(x)^{3}dx}
\end{equation}
converges in measure, as $N\rightarrow\infty$, with respect to the measure $\theta$.
Denote by $G(u)$ the limit of (\ref{cinm}) as $N\rightarrow\infty$.
Then for every $p\in [1,\infty[$,
$G(u)\in L^{p}(d\theta(u))$ and if we set
$
d\rho(u)\equiv G(u)d\theta(u)
$
then the sequence $d\rho_{N}$ converges weakly to $d\rho$ as $N$ tends to
infinity. More precisely for every continuous bounded function
$h:H^{-\sigma}_{0}(S^1)\rightarrow\R$
one has
$$
\int_{H^{-\sigma}_{0}(S^1)}h(u)d\rho(u)=\lim_{N\rightarrow\infty}
\int_{H^{-\sigma}_{0}(S^1)}h(u)d\rho_{N}(u)\,.
$$
\end{theoreme}
Our approach to establish Theorem~\ref{thm} is inspired by the
considerations in \cite{Bo2}.
The main point in the proof of Theorem~\ref{thm} is that thanks to the mean value conservation for
(\ref{BO_N}) the resonant part of $\int_{S^1}u_{N}^{3}$ disappears and thus we
can get the needed integrability by using some known estimates of the second and third
order Wiener chaos. Observe that in a similar analysis in the context of the
2D NLS \cite{Bo2}, the resonant part of the Hamiltonian should be subtracted
which leads to a change of the power nonlinearity to a nonlocal one (the Wick ordering).
\\

In order to prove that the measure $\rho$, constructed in Theorem~\ref{thm} is indeed an invariant measure for the
Benjamin-Ono equation a significant PDE problem should be resolved.
It would be necessary to establish a well-defined
dynamics of (\ref{BO}) for a typical element on the statistical ensemble.
More precisely, one needs to solve almost surely in $\omega$ the Cauchy
problem of (\ref{BO}) with data (\ref{stat}).
Unfortunately, one can prove that the $L^2(S^1)$ of (\ref{stat}) is
a.s. infinity and thus the $L^2$ well-posedness result of Molinet does not apply for 
this data. However, the expression (\ref{stat}) merely misses to belong to $L^2$
(it belongs a.s. to all $H^s(S^1)$, $s<0$). Recall that a somehow similar situation occurred
in \cite{Bo2} and therefore it is not excluded to construct the flow of (\ref{BO})
with data (\ref{stat}) a.s. in $\omega$. Observe that local existence would
suffice since one may exploit the measure invariance of $\mu_N$ under the flow
of (\ref{BO_N}) to get a.s. global solutions (see \cite{Bo}). 
In the final section of this paper we give several estimates confirming that
one may expect to construct the flow of (\ref{BO}) a.s. for data of type (\ref{stat}).
\\

Let us observe that one can use the ideas of this paper to perform similar
constructions with the higher order conservation lows of the Benjamin-Ono equation
in combination with Molinet's well-posedness analysis. We believe that this provides
invariant measures for the Benjamin-Ono equation living on regular spaces.
One however still needs to use the Tao's gauge transform for the truncated ODE
in order to get uniform continuity properties of the flow map. We plan to
pursue these issues elsewhere.
\\

The remaining part of this paper is organized as follows. In the next section, we
prove several elementary inequalities.
In  Section~3, we recall the hypercontractivity properties of the
Ornstein-Uhlenbeck semi-group.
In Section~4 we prove Theorem~\ref{thm}. In the last section we prove several PDE estimates
related to the random series $\varphi(\omega,x)$ which indicate that one may
conjecture that the flow of the Benjamin-Ono equation may be defined for a
typical element of the statistical ensemble.

\section{Elementary calculus inequalities}
In this section, we collect several calculus inequalities, useful for the
sequel. Similar inequalities were used systematically by many authors in the
context of well-posedness for dispersive equations starting from the work of Kenig-Ponce-Vega \cite{KPV}.
\begin{lemme}\label{elem1}
For every $\varepsilon>0$ there exists $C_{\varepsilon}\in\R$ such that for
every $n\in \Z$,
$$
\sum_{n_1\in\Z\backslash\{0,n\}}\frac{1}{|n_1||n-n_1|}\leq \frac{C_{\varepsilon}}{(1+|n|)^{1-\varepsilon}}\,.
$$
\end{lemme}
\begin{proof}
From the triangle inequality, $|n|\leq |n_1|+|n-n_1|$. Therefore, either
$|n|\leq 2|n_1|$ or $|n|\leq 2|n-n_1|$. Thus it suffices to show that for
every $\varepsilon$ there exists $C_{\varepsilon}\in\R$ such that uniformly in
$n$,
$$
\sum_{n_1\in\Z\backslash\{0,n\}}\frac{1}{|n_1|^{\varepsilon}|n-n_1|}\leq
C_{\varepsilon},\quad \sum_{n_1\in\Z\backslash\{0,n\}}\frac{1}{|n_1||n-n_1|^{\varepsilon}}\leq C_{\varepsilon}.
$$
By a change of the summation $n-n_1\rightarrow m$ we observe that the two
inequalities we have to establish are equivalent. Let us prove the second
one. We consider two cases.
\\
{\bf Case 1.} Consider the summation over $n_1$ such that $|n-n_1|\geq
\frac{1}{2}|n_1|$. Denote by $I$ the contribution of this region to the
summation. Then 
$$
I\leq \sum_{n_1\neq 0}\frac{2^{\varepsilon}}{|n_1|^{1+\varepsilon}}=C_{\varepsilon}<\infty\,.
$$
{\bf Case 2.} Consider the summation over $n_1$ such that $|n-n_1|\leq
\frac{1}{2}|n_1|$. Denote by $II$ the contribution of this region to the
sum. The restriction  $|n-n_1|\leq
\frac{1}{2}|n_1|$ implies that $\frac{2}{3}|n|\leq |n_1|\leq 2|n|$. Thus
$$
II\leq C
\sum_{\stackrel{\frac{2}{3}|n|\leq |n_1|\leq 2|n|}{ n_1\neq 0}}\frac{1}{|n_1|}\leq C\log 3\,.
$$
This completes the proof of Lemma~\ref{elem1}.
\end{proof}
\begin{lemme}\label{elem2}
Let us fix $\varepsilon\in]0,1/4[$. Then there exists $C_{\varepsilon}>0$ such
that for every $\alpha\in\Z$,
$$
\sum_{n\in\Z\backslash\{0,\alpha\}}\frac{1}{|n|^{\frac{3}{2}-\varepsilon} |n-\alpha|^{\frac{1}{2}-\varepsilon}}
\leq \frac{C_{\varepsilon}}{(1+|\alpha|)^{\frac{1}{2}-\varepsilon}}\,.
$$ 
\end{lemme}
\begin{proof}
We can suppose that $\alpha\neq 0$. If $|n-\alpha|\geq \frac{|\alpha|}{2}$
then the contribution of these values of $n$ is bounded by 
$$
\Big(\frac{2}{|\alpha|}\Big)^{\frac{1}{2}-\varepsilon}\sum_{n\neq
  0}\frac{1}{|n|^{\frac{3}{2}-\varepsilon}}\leq
\frac{C_{\varepsilon}}{|\alpha|^{\frac{1}{2}-\varepsilon}}\,.
$$
Let us next bound the contribution of those $n$ satisfying $|n-\alpha|\leq
\frac{|\alpha|}{2}$. In this case
$
\frac{|\alpha|}{2}\leq |n|\leq \frac{3|\alpha|}{2}
$
and the contribution of those $n$ to the sum is bounded by
$$
\frac{C_{\varepsilon}}{|\alpha|^{\frac{1}{2}-\varepsilon}}
\sum_{\frac{|\alpha|}{2}\leq |n|\leq \frac{3|\alpha|}{2}  }\frac{1}{|n| |n-\alpha|^{\frac{1}{2}-\varepsilon}}
\leq
\frac{C_{\varepsilon}}{|\alpha|^{\frac{1}{2}-\varepsilon}} 
\sum_{ \frac{|\alpha|}{2}\leq |n|\leq \frac{3|\alpha|}{2}   }\frac{1}{|n|}
\leq
\frac{C_{\varepsilon}}{|\alpha|^{\frac{1}{2}-\varepsilon}}
$$
exactly as in the proof of Lemma~\ref{elem1}. This completes the proof of
Lemma~\ref{elem2}.
\end{proof}
\section{Hypercontractivity properties of the Ornstein-Uhlenbeck semi-group}
In this section, we review some $L^p-L^q$ estimates for the heat flow associated
to the Hartree-Fock operator $\Delta-x\cdot\nabla$ (see Proposition~\ref{hyper} below). We
then obtain corollaries, known as bounds on the Wiener chaos, useful for the
proof of Theorem~\ref{thm}.
For details and background concerning the discussion of this section
(in particular for the proof of  Proposition~\ref{hyper}), we refer
to \cite{ABCFGMRS,LT} and the references therein.
\\

For $d\geq 1$ an integer, we consider the Hilbert space
$H\equiv L^{2}(\R^d,\exp(-|x|^2/2)dx)$ of functions on the euclidean space, square integrable
with respect to the Gaussian measure. Then the operator
\begin{equation}\label{L} 
L\equiv \Delta-x\cdot\nabla=\sum_{j=1}^{d}\Big(\frac{\partial^2}{\partial x_j^{2}}-x_{j}\frac{\partial}{\partial x_j}\Big)
\end{equation}
can be defined as the self adjoint realisation on $L^{2}(\R^d,\exp(-|x|^2/2)dx)$ of $\Delta-x\cdot\nabla$
with domain
$$
D\equiv \Big(u\, :\, u(x)=e^{|x|^2/4}v(x), \quad v\in D_1\Big),
$$
where
$$
D_1\equiv\Big(v\in L^2(\R^d)\,:\,x^{\alpha}\partial^{\beta}v(x)\in L^2(\R^d),
\quad \forall (\alpha,\beta)\in\N^{2d}, |\alpha|+|\beta|\leq 2\Big).
$$
Indeed, one can directly check that
\begin{equation}\label{conjug}
e^{-|x|^2/4}\,L\, e^{|x|^2/4}=\Delta-\Big(\frac{|x|^2}{4}-\frac{d}{2}\Big).
\end{equation}
Of course, $x^2/4$ should be seen as $1/2\int^x ydy$. 
It is well known that $\Delta-|x|^2$ with domain $C^{\infty}_{0}(\R^d)$ is
essentially self adjoint on $L^2(\R^d)$. Moreover $D_1$ defined above is the domain of the
self adjoint extension. In addition 
$$
{\rm spec}(\Delta-|x|^2)=\Big\{-\sum_{j=1}^{d}(2k_j+1),\quad
k_j\in\{0,1,2\cdots\},\quad j=1,\cdots, d \Big\}\,.
$$
We now observe that if $u$ solve $(\Delta-|x|^2)u=\lambda u$ then  $v(x)\equiv
u(x/\sqrt{2})$ solves
$$
\Big(\Delta-\big(\frac{|x|^2}{4}-\frac{d}{2}\big)\Big)v=\frac{\lambda+d}{2}v\,.
$$
Thus, we deduce that
$\Delta-(|x|^2/4-d/2)$, with domain $D_1$, is self-adjoint on $L^2(\R^d)$ and
its spectrum is formed by the
integers $\leq 0$. Therefore, using (\ref{conjug}), we obtain that $L$ has a  self adjoint
realisation on $L^{2}(\R^d,\exp(-|x|^2/2)dx)$ with domain $D$.
The operator $L$ is negative with respect to the
$L^{2}(\R^d,\exp(-|x|^2/2)dx)$ scalar product.
Then the solutions of the linear PDE 
\begin{equation}\label{heat}
\partial_{t}u=L u,\quad u|_{t=0}=u_{0}(x)\in H,\quad x\in\R^d,\,\,
t\in \R^{+}
\end{equation}
are given by the functional calculus of self adjoint operators by the
semi-group $S(t)=\exp(tL)$, i.e. the solution of (\ref{heat}) is given by
$u(t)=S(t)u_0$. Of course one may also define $S(t)$ via the Hille-Yosida theorem.
It turns out that $S(t)$ satisfies an amazing ``smoothing'' property in the
scale of $L^{p}(\R^d,d\mu_d)$, $p\geq 2$, where
$$
d\mu_{d}(x)=(2\pi)^{-d/2}\exp(-|x|^2/2)dx
$$ 
(a probability measure on $\R^d$).
More precisely, a solution starting
from $L^{2}(\R^d,d\mu_d)$ initial data belongs to any $L^{p}(\R^d,d\mu_d)$, $p>2$ (a space smaller than $H$)
for sufficiently long times.
Here is the precise statement. 
\begin{proposition}\label{hyper}
Let us fix $p\geq 2$. Then for every $u_0\in H$, every $t$ satisfying $t\geq \frac{1}{2}\log(p-1)$,
\begin{equation}\label{apriori}
\|S(t)u_0\|_{L^p(\R^d,d\mu_d)}\leq \|u_0\|_{L^2(\R^d,d\mu_d)}\,.
\end{equation}
\end{proposition}
\begin{remarque}
The exponent $2$ in the right hand-side of (\ref{apriori}) may be substituted
by other values $q<p$ and then the restriction on $t$ is $t\geq
(1/2)\log((p-1)/(q-1))$.
There is a close correspondence between (\ref{apriori}) and
logarithmic Sobolev inequalities for the Gaussian measure. In addition,
hypercontractivity estimates of the spirit of (\ref{apriori}) are known for
many other heat flows.  
\end{remarque}

Thanks to (\ref{conjug}) the spectrum of $L$ is formed by the
integers $\leq 0$ and the eigenfunctions of
$L$ may be described in terms of the Hermite polynomials. The Hermite
polynomial $h_{k}(x)$, $k=0,1,2,\cdots$ can be defined via a generating
function as
$$
\exp\big(-\lambda x-\frac{\lambda^2}{2}\big)=\sum_{k=0}^{\infty}\frac{\lambda^k}{\sqrt{k!}}h_{k}(x)\,.
$$
Notice that $h_{0}(x)=1$, $h_{1}(x)=-x$, $h_{2}(x)=\frac{1}{\sqrt{2}}(x^2-1)$. In what follows, we will only need
these three facts about the Hermite polynomials.
A bases of eigenfunctions of $L$ on $H$ is given by 
$$
{\bf h_k}(x)=h_{k_1}(x_1)h_{k_2}(x_2)\cdots h_{k_d}(x_d),
$$
where
$
{\bf k}=(k_1,k_2,\cdots, k_d)\in \N^d
$
and
$ x=(x_1,x_2,\cdots, x_d)\in\R^d\,.$
The eigenfunction ${\bf h_k}$ corresponds to the eigenvalue 
$$
\lambda_{{\bf k}}=-(k_1+\cdots+ k_{d})\,.
$$
The following statement will be used in the proof of
Theorem~\ref{thm}.
\begin{proposition}\label{stein1}
Set
$$
\Sigma_d\equiv 
\big((n_1,n_2,n_3)\in\{1,\cdots, d\}^{3}\, :\,
n_1\neq n_2,n_1\neq n_3,n_2\neq n_3\big).
$$
Then
$$
\|
H(x)
\|_{L^p(\R^d,d\mu_d)}
\leq (p-1)^{\frac{3}{2}}
\|
H(x)
\|_{L^2(\R^d,d\mu_d)},
$$
where
$$
H(x)=\sum_{(n_1,n_2,n_3)\in \Sigma_{d}}
c(n_1,n_2,n_3)x_{n_1}x_{n_2}x_{n_3},\quad c(n_1,n_2,n_3)\in\R\,.
$$
\end{proposition}
\begin{proof}
The function $H$ is an eigenfunction of $L$ corresponding to an eigenvalue
$-3$. Therefore
$
S(t)H=e^{-3t}H.
$
Thus  Proposition~\ref{hyper} yields the bound
$$
\|H\|_{L^p(\R^d,d\mu_d)}\leq
\exp(3t)
\|H\|_{L^2(\R^d,d\mu_d)}\,.
$$
provided $t\geq \frac{1}{2}\log(p-1)$.
By taking $t= \frac{1}{2}\log(p-1)$ in the above bound, we complete the proof of Proposition~\ref{stein1}.
\end{proof}
Let us state another bound related to third order Wiener chaos.
\begin{proposition}\label{stein1.5}
Set
$$
\Sigma_d\equiv 
\big((n_1,n_2)\in\{1,\cdots, d\}^{2}\, :\,
n_1\neq n_2\big).
$$
Then
$$
\|H(x)\|_{L^p(\R^d,d\mu_d)}\leq (p-1)^{\frac{3}{2}}\|H(x)\|_{L^2(\R^d,d\mu_d)},
$$
where
$$
H(x)=\sum_{(n_1,n_2)\in \Sigma_{d}}
c(n_1,n_2)x_{n_1}(x_{n_2}^2-1),\quad c(n_1,n_2)\in\R\,.
$$
\end{proposition}
\begin{proof}
Again the function $H$ is an eigenfunction of $L$ corresponding to an eigenvalue
$-3$. Therefore we can complete the proof as we did in  the proof of Proposition~\ref{stein1}.
\end{proof}
We will also make use of the following inequality.
\begin{proposition}\label{stein2}
We have the bound
$$
\|H(x)\|_{L^p(\R^d,d\mu_d)}\leq (p-1)\|H(x)\|_{L^2(\R^d,d\mu_d)},
$$
where
$$
H(x)=\sum_{n=1}^{d}
c(n)(x_{n}^{2}-1),\quad c(n)\in\R\,.
$$
\end{proposition}
\begin{proof}
The function $H$ is an eigenfunction of $L$ corresponding to an eigenvalue
$-2$. Therefore
$
S(t)H=e^{-2t}H.
$
Thus  Proposition~\ref{hyper} yields the bound
$$
\|H\|_{L^p(\R^d,d\mu_d)}\leq
\exp(2t)
\|H\|_{L^2(\R^d,d\mu_d)}\,.
$$
provided $t\geq \frac{1}{2}\log(p-1)$.
As in the proof of  Proposition~\ref{stein1}
by taking $t= \frac{1}{2}\log(p-1)$ in the above bound, we complete the proof.
\end{proof}
\section{Proof of  Theorem~\ref{thm}}
In order to deal with the low frequencies we will need the following
distributional inequality.
\begin{proposition}\label{calais}
For every $C_1>0$ and $C_2>0$, $\varepsilon>0$, $\alpha>0$ there exist
$C>0$, $c>0$ such that for every integer $N\geq 1$, every $\lambda\geq 2$ satisfying $N\leq\lambda^{\alpha}$ one has
$$
\theta\Big(u\in H^{-\sigma}_{0}(S^1)\,:\, \|S_{N}u\|_{L^{\infty}(S^1)}\geq
C_1\lambda,\, \|S_{N}u\|^{2}_{L^2(S^1)}\leq C_2\log\lambda\Big)\leq \frac{C}{\exp(c\lambda^{2-\varepsilon})}\,.
$$
\end{proposition}
\begin{proof}
We will need the following Khinchin type inequality.
\begin{lemme}\label{lem1}
Let $(l_n(\omega))_{n\in\N}$ be a sequence of independent identically
distributed standard real
Gaussian random variables. Then for every
$\lambda >0$, every sequence $(c_n)\in l^2(\N)$
of real numbers,
\begin{equation}\label{khin}
p\Big(\omega\,:\,\big|\sum_{n=0}^{\infty} c_n l_n(\omega) \big|>\lambda\Big)\leq
2 e^{-\frac{\lambda^2}{2\sum_{n}c_n^2}}
\, .
\end{equation}
\end{lemme}
\begin{proof}
The assertion of this lemma follows from the estimates on first order Wiener
chaos considered in the previous section. 
It is also a consequence of the observation that $\sum c_n l_n$ is a Gaussian
in ${\mathcal N}(0,\sigma^2)$ with $\sigma^2=\sum_{n}c_n^2$.
We include however here a proof of (\ref{khin}) which has the
advantage to work for more general systems of independent zero mean
value random variables instead of $(l_n(\omega))_{n\in\N}$ (such as Bernouli variables).

For $t>0$ to be determined later, using the independence, we obtain that
\begin{eqnarray*}
\int_{\Omega}\, e^{t\sum_{n\geq 0}c_n l_n(\omega)}dp(\omega)
& = &\prod_{n\geq 0}\int_{\Omega}e^{t c_n l_n(\omega)}dp(\omega)
\\
& = & \prod_{n\geq 0}\int_{-\infty}^{\infty}e^{tc_n x}\, e^{-x^2/2}\frac{dx}{\sqrt{2\pi}}
\\
& = &
\prod_{n\geq 0}e^{(t c_n)^2/2}= e^{(t^2/2)\sum_{n}c_n^2}\, .
\end{eqnarray*}
Using the above calculation, we infer that
$$
e^{(t^2/2)\sum_{n}c_n^2}\geq e^{t\lambda}\,\,\,  p\,(\omega\,:\,\sum_{n\geq 1} c_n l_n(\omega)>\lambda)
$$
or equivalently,
$$
p\,(\omega\,:\,\sum_{n\geq 1} c_n l_n(\omega)>\lambda)\leq e^{(t^2/2)\sum_{n}c_n^2}\,\,\,
e^{-t\lambda}\, .
$$
Using that for $a>0$ the minimum of $f(t)=at^2-bt$ is $-b^2/4a$, we obtain
that 
$$
p\,(\omega\,:\,\sum_{n\geq 1} c_n l_n(\omega)>\lambda)\leq
e^{-\frac{\lambda^2}{2\sum_{n}c_n^2}}\, .
$$
In the same way (replacing $c_n$ by $-c_n$), we can show that
$$
p\,(\omega\,:\,\sum_{n\geq 1} c_n l_n(\omega)<-\lambda)\leq
e^{-\frac{\lambda^2}{2\sum_{n}c_n^2}}
$$
which completes the proof of Lemma~\ref{lem1}.
\end{proof}
Let us now give the proof of Proposition~\ref{calais}. Set
$$
A_{\lambda}\equiv (u\in H^{-\sigma}_{0}(S^1)\,:\, \|S_{N}u\|_{L^{\infty}(S^1)}\geq
C_1\lambda,\, \|S_{N}u\|^{2}_{L^2(S^1)}\leq C_2\log\lambda)\,.
$$
Observe that for $\alpha<2$ the Sobolev embedding applied to $S_N u$ suffices to
conclude that for $\lambda\gg 1$ the set $A_{\lambda}$ is empty. Hence the
result is not trivial for $\alpha\geq 2 $ (which will be the case in our
application of Proposition~\ref{calais}).
For $\alpha\geq 2 $ the  Sobolev embedding applied to $S_N u$ does not give a
lower bound for $\|S_N u\|_{L^2}$ which is the main source of difficulty. 
Let us fix $\beta>2\alpha$.
Define the points $x_j\in S^1$, $j=0,\cdots, [\Lambda\lambda^\beta]$, where
$\Lambda\gg 1$ is to be fixed later by
$
x_{j}\equiv (2\pi j)/(\Lambda\lambda^{\beta}).
$
The number $\Lambda$ may depend on $C_1$, $C_2$, $\varepsilon$, $\alpha$ but
should be independent of $\lambda$ and $N$.
Notice that ${\rm dist}(x_{j},x_{j+1})\leq 2\pi/(\Lambda\lambda^{\beta})$,
where $x_{[\Lambda\lambda^{\beta}]+1}\equiv x_{0}$ and ${\rm dist}$ denotes the
distance on $S^1$ (i.e. ${\rm mod }\,2\pi$). Next, we define the sets
$A_{\lambda,j}$ by
$$
A_{\lambda,j}\equiv
 (u\in H^{-\sigma}_{0}(S^1)\,:\, |S_{N}u(x_j)|\geq
\frac{1}{2}C_1\lambda,\, \|S_{N}u\|^{2}_{L^2(S^1)}\leq C_2\log\lambda)\,.
$$
We claim that for $\Lambda\gg 1$,
\begin{equation}\label{bruge1}
A_{\lambda}\subset
\bigcup_{j=0}^{[\Lambda\lambda^\beta]}A_{\lambda,j}\,.
\end{equation}
Let us prove (\ref{bruge1}).
Fix $u\in A_{\lambda}$. Let $x^{\star}\in S^1$ be such that
$$
|S_{N}u(x^{\star})|=\max_{x\in S^1}|S_{N}u(x)|.
$$
Thus 
$
|S_{N}u(x^{\star})|\geq C_1\lambda.
$
Then there exists $j_0\in \{0,\cdots, [\Lambda\lambda^\beta]\}$ such that 
$$
|x^\star-x_{j_0}|\leq \frac{2\pi}{\Lambda\lambda^{\beta}}.
$$
Then we can write
\begin{eqnarray*}
|S_{N}u(x^\star)-S_{N}u(x_{j_0})|& = & \Big|\int_{x_{j_0}}^{x^\star}(S_N
u)'(tx^\star+(1-t)x_{j_0})dt\Big|
\\
& \leq &
|x^\star-x_{j_0}|^{\frac{1}{2}}\|(S_N u)'\|_{L^2(S^1)}
\\
& \leq &
\frac{\sqrt{2\pi}}{\sqrt{\Lambda\lambda^{\beta}}}N\|S_N u\|_{L^2(S^1)}
\\
& \leq &
\frac{\sqrt{2\pi}}{\sqrt{\Lambda\lambda^{\beta}}}\lambda^{\alpha}\sqrt{C_2 \log\lambda}\,.
\end{eqnarray*}
Let us choose $\Lambda\gg 1$ such that for every $\lambda\geq 2$,
$$
\frac{\sqrt{2\pi}}{\sqrt{\Lambda\lambda^{\beta}}}\lambda^{\alpha}\sqrt{C_2 \log\lambda}
\leq 
\frac{1}{2}C_1\lambda.
$$
Then by the triangle inequality
$$
|S_{N}u(x_{j_0})|\geq
|S_{N}u(x^\star)|
-|S_{N}u(x^\star)-S_{N}u(x_{j_0})|\geq C_{1}\lambda-\frac{1}{2}C_1\lambda=
\frac{1}{2}C_1\lambda\,.
$$
Hence $u\in A_{\lambda,j_0}$ which proves (\ref{bruge1}). Let us next evaluate
$\theta(A_{\lambda,j})$. For that purpose we will make appeal to Lemma~\ref{lem1}.
Observe that
\begin{multline*}
\theta(A_{\lambda,j})=
p
\Big(\omega\,:\,
\Big|
\sum_{n=1}^{N}(2\pi n)^{-\frac{1}{2}}
\big(
\cos(nx_j)h_{n}(\omega)+\sin(nx_j)l_{n}(\omega)
\big)
\Big|\geq \frac{1}{2} C_1\lambda,
\\
\sum_{n=1}^{N}n^{-1}(h_n^{2}(\omega)+l_{n}^{2}(\omega))\leq 2C_2\log\lambda
\Big).
\end{multline*}
Therefore, by ignoring the $L^2$ restriction and using  Lemma~\ref{lem1}, we
obtain that
$$
\theta(A_{\lambda,j})\leq 2 e^{-\frac{(C_1\lambda)^2}{8\kappa}}\,,
$$
where
$$
\kappa=\sum_{n=1}^{N}\Big(
\frac{\cos^2(nx_j)}{2\pi n}+\frac{\sin^2(nx_j)}{2\pi n}
\Big)=\frac{1}{2\pi}\sum_{n=1}^{N}\frac{1}{n}\,.
$$
Thus using that $N\leq \lambda^{\alpha}$, we infer that
$\kappa \leq C\log\lambda$, where $C$ is independent of $N$ and $\lambda$.
Therefore there exists $c>0$, depending only on $C_1,C_2,\alpha,\varepsilon$, such
that 
\begin{equation}\label{bruge2}
\theta(A_{\lambda,j})\leq 2e^{-c\lambda^{2-\varepsilon/2}}\,.
\end{equation}
Combining (\ref{bruge1}) and (\ref{bruge2}) implies that
$$
\theta(A_{\lambda})\leq
\sum_{j=0}^{[\Lambda\lambda^\beta]}
\theta(A_{\lambda,j})\leq
2(\Lambda\lambda^\beta+1) e^{-c\lambda^{2-\varepsilon/2}}
\leq C e^{-c\lambda^{2-\varepsilon}}\,,
$$
where $C,c>0$ are independent of $\lambda$ and $N$.
This completes the proof of Proposition~\ref{calais}.
\end{proof}
Let us define the functions
$f_{N}\,:\,H^{-\sigma}_{0}(S^1)\rightarrow\R$
by
$$
f_{N}(u)\equiv \int_{S^1}((S_{N}u)(x))^{3}dx\,.
$$
Then we have the following statement. 
\begin{lemme}\label{integral}
The sequence $(f_N)_{N\geq 1}$ is a Cauchy sequence in
$L^2(H^{-\sigma}_{0}(S^1),{\mathcal B},d\theta)$.
More precisely, for every $\alpha<1/2$ there exists $C>0$ such that for every $M>N\geq 1$,
\begin{equation}\label{bound}
\Big\|
f_{M}(u)-f_{N}(u)
\Big\|_{L^2(H^{-\sigma}_{0}(S^1),{\mathcal B},d\theta)}
\leq
CN^{-\alpha}\,.
\end{equation}
Moreover, for every $M>N\geq 1$,
every $p\geq 2$,
\begin{equation}\label{boundpak}
\Big\|
f_{M}(u)-f_{N}(u)
\Big\|_{L^p(H^{-\sigma}_{0}(S^1),{\mathcal B},d\theta)}
\leq
Cp^{\frac{3}{2}}N^{-\alpha}\,.
\end{equation}
\end{lemme}
Denote by $f(u)\in  L^2(H^{-\sigma}_{0}(S^1),{\mathcal B},d\theta)$ the limit
of $(f_{N})_{N\geq 1}$. 
Let us notice that the result of Lemma~\ref{integral} is 
displaying some important cancellations
since using (for instance) the Fernique integrability theorem one may show 
that $\int_{S^1}|u|^{3}=\infty$, $\theta$ a.s.
\begin{proof}[Proof of Lemma~\ref{integral}]
Write
\begin{eqnarray*}
\|f_{N}\|^{2}_{L^2(H^{-\sigma}_{0}(S^1),{\mathcal B},d\theta)}
& = &
\int_{H^{-\sigma}_{0}(S^1)}
\Big|
\int_{S^1}((S_{N}u)(x))^{3}dx
\Big|^{2}d\theta(u)
\\
& = &
\int_{\Omega}
\Big|
\int_{S^1}((S_{N}\varphi(\omega,x)))^{3}dx
\Big|^{2}dp(\omega),
\end{eqnarray*}
where $\varphi(\omega,x)$ is defined by (\ref{stat}).
For $N\geq 2$, we set
$$
\Sigma(N)\equiv
\big(
(n_1,n_2,n_3)\in\Z^3\, :\,
n_1+n_2+n_3=0,\,
0<|n_1|,|n_2|,|n_3|\leq N
\big).
$$
Then
$$
\int_{S^1}(S_{N}\varphi(\omega,x))^3dx=
\frac{1}{8\pi^{\frac{3}{2}}}\sum_{(n_1,n_2,n_3)\in \Sigma(N)}
\frac{g_{n_1}(\omega)}{\sqrt{|n_1|}}\frac{g_{n_2}(\omega)}{\sqrt{|n_2|}}\frac{g_{n_3}(\omega)}{\sqrt{|n_3|}}\,.
$$
Next we define $\Sigma_{1}(N)$ as follows
$$
\Sigma_{1}(N)\equiv 
\big((n_1,n_2,n_3)\in\Sigma(N)\, :\,
n_1\neq\pm n_2,n_1\neq \pm n_3,n_2\neq\pm n_3\big).
$$
Observe that triples of the form $(n,-n,0)$ can not belong to $\Sigma(N)$
and therefore, we may write
$$
\int_{S^1}(S_{N}\varphi(\omega,x))^3=F_{1}(N,\omega)+F_{2}(N,\omega),
$$
$$
F_{1}(N,\omega)\equiv \frac{3}{8\pi^{\frac{3}{2}}} 
\sum_{0<|n|\leq N/2}\frac{g_{n}^{2}(\omega)\overline{g_{2n}(\omega)}}{|n|^{\frac{3}{2}}}
$$
is the contribution of the terms $(n,n,-2n)$, $(n,-2n,n)$ and $(-2n,n,n)$ and
$$
F_{2}(N,\omega)\equiv 
\frac{1}{8\pi^{\frac{3}{2}}}\sum_{(n_1,n_2,n_3)\in \Sigma_1(N)}
\frac{g_{n_1}(\omega)}{\sqrt{|n_1|}}\frac{g_{n_2}(\omega)}{\sqrt{|n_2|}}\frac{g_{n_3}(\omega)}{\sqrt{|n_3|}}\,.
$$
is the contribution of the remaining terms.
Since
\begin{multline*}
\|f_{M}-f_{N}\|^{2}_{L^2(H^{-\sigma}_{0}(S^1),{\mathcal B},d\theta)}
\\
=
\int_{\Omega}
\Big|
\int_{S^1}((S_{M}\varphi(\omega,x)))^{3}dx
-
\int_{S^1}((S_{N}\varphi(\omega,x)))^{3}dx
\Big|^{2}dp(\omega),
\end{multline*}
it suffices to show that $(F_j(N,\cdot))_{N\geq 1}$, $j=1,2$ are Cauchy
sequences in $L^2(\Omega)$ satisfying bounds of type (\ref{bound}), (\ref{boundpak}).
Using the H\"older inequality in the $\Omega$ integration, we may write
\begin{eqnarray*}
\|F_{1}(M,\omega)-F_{1}(N,\omega)\|_{L^2(\Omega)}
& \leq & 
C
\sum_{N/2<|n|\leq M/2}\frac{\|g_{n}\|_{L^3(\Omega)}^{2}\|g_{2n}\|_{L^3(\Omega)}}
{|n|^{\frac{3}{2}}}
\\
& \leq & C\sum_{N/2<|n|\leq M/2}\frac{1}{|n|^{\frac{3}{2}}}
\\
& \leq & 
\frac{C_{\alpha}}{N^{\alpha}}\, ({\rm recall\,\, that} \,\, \alpha<1/2) \,.
\end{eqnarray*}
Thus $(F_1(N,\cdot))_{N\geq 1}$ is a Cauchy sequence in $L^2(\Omega)$ with the
needed quantitative bound.
Let us next analyse $F_{2}(N,\omega)$.
For that purpose, in contrast with $F_{1}(N,\omega)$, an orthogonality argument
will be needed.
For $M>N\geq 1$, we set
\begin{multline*}
\Lambda(N,M)\equiv
\big(
(n_1,n_2,n_3)\in\Z^3\, :\,
n_1+n_2+n_3=0,\, n_1\neq\pm n_2,n_1\neq \pm n_3,n_2\neq\pm n_3
\\ 
0<|n_1|,|n_2|,|n_3|\leq M,\,
\max(|n_1|,|n_2|,|n_3|)>N
\big).
\end{multline*}
Therefore, we can write
$$
F_{2}(M,\omega)-F_{2}(N,\omega)
=
\frac{1}{8\pi^{\frac{3}{2}}}\sum_{(n_1,n_2,n_3)\in \Lambda(N,M)}
\frac{g_{n_1}(\omega)}{\sqrt{|n_1|}}\frac{g_{n_2}(\omega)}{\sqrt{|n_2|}}\frac{g_{n_3}(\omega)}{\sqrt{|n_3|}}\,.
$$
Observe that if $(n_1,n_2,n_3)$ and $(m_1,m_2,m_3)$ are two triples
from $\Lambda(N,M)$ such that $\{n_1,n_2,n_3\}\neq \{m_1,m_2,m_3\}$ then
\begin{equation}\label{triples}
\int_{\Omega}g_{n_1}(\omega)g_{n_2}(\omega)g_{n_3}(\omega)\overline{g_{m_1}(\omega)g_{m_2}(\omega)g_{m_3}(\omega)}dp(\omega)=0.
\end{equation}
Indeed, using the independence, if $n_{j_1}=-m_{j_2}$ for some
$j_1,j_2\in\{1,2,3\}$ then the integral (\ref{triples}) is zero since
$\int_{\Omega}g^2_{n_{j_1}}(\omega)dp(\omega)=0$
and $\pm n_{j_1}$ can not belong to the remaining indexes. In all other cases
there is one of the indexes $(n_1,n_2,n_3,m_1,m_2,m_3)$ which is repeated
only once and its opposite does not belong to $\{n_1,n_2,n_3,m_1,m_2,m_3\}$.
Therefore, we can write
\begin{eqnarray*}
\|F_{2}(M,\omega)-F_{2}(N,\omega)\|_{L^2(\Omega)}^{2}
& \leq &
C\sum_{(n_1,n_2,n_3)\in \Lambda(N,M)}
\frac{1}{|n_1|}
\frac{1}{|n_2|}
\frac{1}{|n_3|}
\\
& \leq &
C\sum_{n_1\in \Z}\sum_{|n_2|\geq N}
\frac{1}{(1+|n_1|)(1+|n_2|)(1+|n_1+n_2|)}\,.
\end{eqnarray*}
Using Lemma~\ref{elem1}, we infer that
$$
\sum_{n_1\in \Z}
\frac{1}{(1+|n_1|)(1+|n_1+n_2|)}
\leq \frac{C_{\varepsilon}}{(1+|n_2|)^{1-\varepsilon}}
+
\frac{2}{(1+|n_2|)}
<
\frac{C_{\varepsilon}+2}{(1+|n_2|)^{1-\varepsilon}}\,.
$$
Therefore
$$
\|
F_{2}(M,\omega)-F_{2}(N,\omega)
\|_{L^2(\Omega)}^{2}\leq C
\sum_{|n_2|\geq N}\frac{1}{(1+|n_2|)^{2-\varepsilon}}\leq \frac{C_{\alpha}}{N^{2\alpha}}\,,
$$
provided $1-2\alpha>\varepsilon>0$.
Therefore $(F_{2}(\omega,N))_{N\geq 1}$ is a Cauchy sequence in $L^2(\Omega)$
with the needed quantitative bound.
This completes the proof of (\ref{bound}).
Let us now turn to the proof of (\ref{boundpak}).
Write via the triangle inequality,
\begin{multline*}
\|f_{M}-f_{N}\|^{p}_{L^p(H^{-\sigma}_{0}(S^1),{\mathcal B},d\theta)}
\\
=
\int_{\Omega}
\Big|
\int_{S^1}((S_{M}\varphi(\omega,x)))^{3}dx
-
\int_{S^1}((S_{N}\varphi(\omega,x)))^{3}dx
\Big|^{p}dp(\omega)
\\
\leq
\Big(
\|F_{1}(M,\omega)-F_{1}(N,\omega)\|_{L^p(\Omega)}
+\|F_{2}(M,\omega)-F_{2}(N,\omega)\|_{L^p(\Omega)}
\Big)^{p}\,.
\end{multline*}
Write
$$
F_{1}(M,\omega)-
F_{1}(N,\omega)
= \frac{3}{4\pi^{\frac{3}{2}}} 
\sum_{N/2<n\leq M/2}
\frac{{\rm Re}\big(g_{n}^{2}(\omega)\overline{g_{2n}(\omega)}\big)}{|n|^{\frac{3}{2}}}\,.
$$
Recall that
$
g_{n}(\omega)=\frac{1}{\sqrt{2}}(h_{n}(\omega)-il_{n}(\omega))
$
and thus one may directly check that
$$
{\rm Re}\big(g_{n}^{2}(\omega)\overline{g_{2n}(\omega)}\big)=
\frac{1}{2\sqrt{2}}
\Big(
(h^2_{n}(\omega)-1)h_{2n}(\omega)-(l_{n}^{2}(\omega)-1)h_{2n}(\omega)
+2h_{n}(\omega)l_{n}(\omega)l_{2n}(\omega)
\Big).
$$
Hence we are in the scope of applicability of Propositions~\ref{stein1} and \ref{stein1.5}.
Consider $\R^{2M}$ parametrized
by $(x_1,\cdots,x_{M},y_1,\cdots,y_{M})$, where $(x_1,\cdots,x_{M})$
correspond to the $h_{n}(\omega)$, $n=1,\cdots,M$ and where $(y_1,\cdots,y_{M})$
correspond to the $l_{n}(\omega)$, $n=1,\cdots,M$. 
Then we will apply Proposition~\ref{stein1.5} (with $d=2M$) to the function
$$
H_1(x_1,\cdots,x_{M},y_1,\cdots,y_{M})
\equiv\frac{3}{8\sqrt{2}\pi^{\frac{3}{2}}}
\sum_{N/2<n\leq M/2}|n|^{-\frac{3}{2}}
\big((x_{n}^2-1)x_{2n}-(y_{n}^2-1)x_{2n}\big)
$$ 
and Proposition~\ref{stein1} to the function
$$
H_2(x_1,\cdots,x_{M},y_1,\cdots,y_{M})
\equiv\frac{3}{4\sqrt{2}\pi^{\frac{3}{2}}}\sum_{N/2<n\leq M/2}
|n|^{-\frac{3}{2}}
x_{n}y_{n}y_{2n}\,.
$$
Indeed,
\begin{eqnarray*}
F_{1}(M,\omega)-F_{1}(N,\omega) & = &
H_1(h_1(\omega),\cdots,h_{M}(\omega),l_1(\omega),\cdots,l_{M}(\omega))
\\
& &
+
H_2(h_1(\omega),\cdots,h_{M}(\omega),l_1(\omega),\cdots,l_{M}(\omega)).
\end{eqnarray*}
Using the independence, we may write that for $j=1,2$,
\begin{multline*}
\|H_j(h_1(\omega),\cdots,h_{M}(\omega),l_1(\omega),\cdots,l_{M}(\omega))\|_{L^p(\Omega)}
=
\\
=
\|H_j(x_1,\cdots,x_{M},y_1,\cdots,y_{M})\|_{L^p\Big(\R^{2M},(2\pi)^{-M}\exp\big(-\frac{1}{2}\sum_{n=1}^{M}(x_n^2+y_n^2)\big)dx_1...dy_M\Big)}\,.
\end{multline*}
Therefore, using Proposition~\ref{stein1.5} and  Proposition~\ref{stein1},
by splitting $F_{1}(M,\omega)- F_{1}(N,\omega)$ into two parts, we obtain that
$$
\|F_{1}(M,\omega)-F_{1}(N,\omega)\|_{L^p(\Omega)}\leq 
Cp^{3/2}\|F_{1}(M,\omega)-F_{1}(N,\omega)\|_{L^2(\Omega)}
\leq C_{\alpha}p^{3/2}N^{-\alpha}\,,
$$
where $C_{\alpha}$ is independent of $p$, $M$ and $N$.
Similarly, by developing the product
$$
g_{n_1}(\omega)g_{n_2}(\omega)g_{n_3}(\omega)
$$
for
$
(n_1,n_2,n_3)\in \Lambda(N,M)
$
we observe that the difference $F_{2}(M,\omega)-F_{2}(N,\omega)$ 
fits in the scope of applicability of Proposition~\ref{stein1}. We obtain that
$$
\|F_{2}(M,\omega)-F_{2}(N,\omega)\|_{L^p(\Omega)}\leq 
Cp^{3/2}\|F_{2}(M,\omega)-F_{2}(N,\omega)\|_{L^2(\Omega)}
\leq C_{\alpha}p^{3/2}N^{-\alpha}\,.
$$
Thus (\ref{boundpak}) is established.
This completes the proof of Lemma~\ref{integral}.
\end{proof}
We have the following standard corollary of  Lemma~\ref{integral}.
\begin{corollaire}\label{cheb}
Under the assumption of Lemma~\ref{integral},
the sequence $(f_{N})_{N\geq 1}$ converges in measure to $f$. More
precisely, for every $\varepsilon>0$,
$$
\lim_{N\rightarrow\infty}
\theta(u\in H^{-\sigma}_{0}(S^1)\,:\,
|f(u)-f_{N}(u)|>\varepsilon)=0.
$$
\end{corollaire}
\begin{proof}
This is a consequence of the Tchebishev inequality.
\end{proof}
The next lemma is a general feature.
\begin{lemme}\label{general}
Let $F$ be a real valued measurable function on $H^{-\sigma}_{0}(S^1)$. Suppose
that there exist $\alpha>0$, $N>0$, $k\in \N^{\star}$ and $C>0$ such that for every $p\geq 2$ one has
\begin{equation}\label{hyp}
\|F\|_{L^p(d\theta)}\leq CN^{-\alpha}p^{k/2}\,.
\end{equation}
Then there exists $\delta>0$ and $C_1>0$ depending on $C$ and $k$
but {\bf independent} of $N$ and $\alpha$ such that
\begin{equation}\label{golse}
\int_{H^{-\sigma}_{0}(S^1)}e^{\delta N^{\frac{2\alpha}{k}}|F(u)|^{\frac{2}{k}}}d\theta(u)\leq C_1.
\end{equation}
As a consequence for $\lambda>0$,
\begin{equation}\label{nicolas}
\theta(u\in H^{-\sigma}_{0}(S^1)\,:\, |F(u)|>\lambda)\leq C_{1}e^{-\delta N^{\frac{2\alpha}{k}}\lambda^{\frac{2}{k}}}\,.
\end{equation}
\end{lemme}
\begin{proof}
If one is only interested to get (\ref{nicolas}) then it suffices to use the
Tchebishev inequality in the context of (\ref{hyp}) with a suitable $p$ (depending of $\lambda$).
Let us now give the proof of the claimed statement (\ref{golse}).
Write
$$
e^{\delta N^{\frac{2\alpha}{k}}|F(u)|^{\frac{2}{k}}}
=
\sum_{n=0}^{k-1}
\frac{\delta^{n}N^{\frac{2\alpha n}{k}}|F(u)|^{\frac{2n}{k}}}{n!}
+
\sum_{n=k}^{\infty}
\frac{\delta^{n}N^{\frac{2\alpha n}{k}}|F(u)|^{\frac{2n}{k}}}{n!}\,.
$$
If $k\geq 2$, using the H\"older inequality and (\ref{hyp}), we get for
$n=1,\cdots,k-1$,
$$
\int_{H^{-\sigma}_{0}(S^1)}|F(u)|^{\frac{2n}{k}}d\theta(u)\leq
\|F\|_{L^{2n}(d\theta)}^{\frac{2n}{k}}\leq
\Big[CN^{-\alpha}(2n)^{\frac{k}{2}}\Big]^{\frac{2n}{k}}
= C^{\frac{2n}{k}}N^{-\frac{2\alpha n}{k}}(2n)^{n}\,.
$$
The Stirling formula provides the existence of a positive constant $\tilde{C}$
such that for every integer $n\geq 1$,
$$
\frac{n^n}{n!}\leq\tilde{C}\frac{e^n}{\sqrt{n}}\,.
$$
Therefore, by using (\ref{hyp}), we obtain that for $n\geq k$,
\begin{eqnarray*}
\int_{H^{-\sigma}_{0}(S^1)}
\frac{\delta^{n}N^{\frac{2\alpha n}{k}}|F(u)|^{\frac{2n}{k}}}{n!}
d\theta(u)
& \leq &
\frac{\delta^{n}N^{\frac{2\alpha n}{k}}}{n!}
\Big[
CN^{-\alpha}\big(\frac{2n}{k}\big)^{\frac{k}{2}}
\Big]^{\frac{2n}{k}}
\\
& = &
\frac{n^n}{n!}
\Big(\frac{2}{k}C^{\frac{2}{k}}\delta\Big)^{n}
\leq
\frac{\tilde{C}}{\sqrt{n}}\Big(\frac{2}{k}C^{\frac{2}{k}}e\delta\Big)^{n}.
\end{eqnarray*}
Summarizing the preceding gives that for $k\geq 2$,
$$
\int_{H^{-\sigma}_{0}(S^1)}
e^{\delta N^{\frac{2\alpha}{k}}|F(u)|^{\frac{2}{k}}}
d\theta(u)
\leq 1+\sum_{n=0}^{k-1}
\frac{C^{\frac{2n}{k}}(2n)^{n}}{n!}\delta^{n}
+
\tilde{C}
\sum_{n=k}^{\infty}\Big(\frac{2}{k}C^{\frac{2}{k}}e\delta\Big)^{n}
\leq C_1,
$$
provided that $\delta>0$ is such that
$$\delta<\frac{k}{2C^{\frac{2}{k}}e}.$$ For $k=1$, the same bound holds by
replacing the term 
$$
\sum_{n=0}^{k-1}
\frac{C^{\frac{2n}{k}}(2n)^{n}}{n!}\delta^{n}
$$
in the above inequality by zero.
This completes the proof of Lemma~\ref{general}.
\end{proof}
Lemma~\ref{general} implies the following distributional inequality for $(f_{N})_{N\geq 1}$.
\begin{lemme}\label{distri}
For every $\alpha<1/2$ there exists $C>0$ and $\delta>0$ such that for every $M>N\geq 1$,
every $\lambda>0$
$$
\theta(u\in H^{-\sigma}_{0}(S^1)\,:\, |f_{M}(u)-f_{N}(u)|>\lambda)
\leq Ce^{-\delta (N^{\alpha}\lambda)^{2/3}}\,.
$$
\end{lemme}
\begin{proof}
It suffices to combine Lemma~\ref{integral} and Lemma~\ref{general}.
\end{proof}
We next study the limit of $\|S_{N}(u)\|_{L^2(S^1)}^{2}-\alpha_{N}$ as
$N\rightarrow\infty$.
Let us define the functions
$g_{N}\,:\,H^{-\sigma}_{0}(S^1)\rightarrow\R$
by
\begin{equation}\label{gN}
g_{N}(u)\equiv \|S_{N}(u)\|_{L^2(S^1)}^{2}-\alpha_{N}\,.
\end{equation}
We have the following statement.
\begin{lemme}\label{Wick}
The sequence $(g_N)_{N\geq 1}$ is a Cauchy sequence in
$L^2(H^{-\sigma}_{0}(S^1),{\mathcal B},d\theta)$.
More precisely, there exists $C>0$ such that for every $M>N\geq 1$,
\begin{equation}\label{boundtris}
\Big\|
g_{M}(u)-g_{N}(u)
\Big\|_{L^2(H^{-\sigma}_{0}(S^1),{\mathcal B},d\theta)}
\leq
CN^{-\frac{1}{2}}\,.
\end{equation}
Moreover, if we denote by $g(u)$ the limit of $g_{N}(u)$ in
$L^2(H^{-\sigma}_{0}(S^1),{\mathcal B},d\theta)$
then $g_{N}(u)$ converges to $g(u)$ in measure :
$$
\forall \varepsilon>0,\quad 
\lim_{N\rightarrow\infty}
\theta(u\in H^{-\sigma}_{0}(S^1)\,:\,
|g(u)-g_{N}(u)|>\varepsilon)=0.
$$
\end{lemme}
\begin{proof}
Write
\begin{multline*}
\|g_{M}-g_{N}\|^{2}_{L^2(H^{-\sigma}_{0}(S^1),{\mathcal B},d\theta)}
\\
=\int_{\Omega}\Big|\|S_{M}\varphi(\omega,\cdot)\|_{L^2(S^1)}^{2}-\alpha_{M}-
\|S_{N}\varphi(\omega,\cdot)\|_{L^2(S^1)}^{2}+\alpha_{N}\Big|^{2}dp(\omega)
\\
=\frac{1}{4}\int_{\Omega}\Big|\sum_{N<|n|\leq M}\frac{|g_{n}(\omega)|^{2}-1}{|n|}\Big|^{2}dp(\omega)\,.
\end{multline*}
Thanks to the independence and the normalization of $(g_{n}(\omega))$ we obtain that for $n_1\neq
n_2$ one has
$$
\int_{\Omega}(|g_{n_1}(\omega)|^{2}-1)(|g_{n_2}(\omega)|^{2}-1)dp(\omega)=0\,.
$$
Therefore
$$
\|g_{M}-g_{N}\|^{2}_{L^2(H^{-\sigma}_{0}(S^1),{\mathcal B},d\theta)}
=
\sum_{N<|n|\leq M}\frac{c}{|n|^2}
\leq \frac{C}{N}\,.
$$
This proves (\ref{boundtris}). The convergence of $(g_{N}(u))$ in measure
follows from the Chebishev inequality.
This completes the proof of Lemma~\ref{Wick}.
\end{proof}
We now prove a distributional inequality for $(g_{N})_{N\geq 1}$.
\begin{lemme}\label{Wickpak}
There exist $C>0$ and $\delta>0$ such that for every $M>N\geq 1$,
every $\lambda>0$
$$
\theta(u\in H^{-\sigma}_{0}(S^1)\,:\, |g_{M}(u)-g_{N}(u)|>\lambda)
\leq Ce^{-\delta N^{\frac{1}{2}}\lambda}\,.
$$
\end{lemme}
\begin{proof}
We have 
$$
\|g_{M}-g_{N}\|^{p}_{L^p(H^{-\sigma}_{0}(S^1),{\mathcal B},d\theta)}=\Big(\frac{1}{2}\Big)^{p}\int_{\Omega}\Big|\sum_{N<|n|\leq
  M}\frac{|g_{n}(\omega)|^{2}-1}{|n|}\Big|^{p}dp(\omega)
$$
Recall that
$
g_{n}(\omega)=\frac{1}{\sqrt{2}}(h_{n}(\omega)-il_{n}(\omega))
$
and thus
$$
|g_{n}(\omega)|^{2}-1=\frac{1}{2}\big(h_{n}^{2}(\omega)-1\big)+\frac{1}{2}\big(l_{n}^{2}(\omega)-1\big).
$$
Therefore, using Proposition~\ref{stein2} and (\ref{boundtris}), we obtain that
$$
\|g_{M}-g_{N}\|_{L^p(H^{-\sigma}_{0}(S^1),{\mathcal B},d\theta)}
\leq Cp\|g_{M}-g_{N}\|_{L^2(H^{-\sigma}_{0}(S^1),{\mathcal B},d\theta)}
\leq 
CpN^{-\frac{1}{2}}\,.
$$
A use of Lemma~\ref{general} completes the proof of Lemma~\ref{Wickpak}.
\end{proof}
Combining Lemma~\ref{integral} and Lemma~\ref{Wick}, we may define the
function 
$$
G\,:\,H^{-\sigma}_{0}(S^1)\longrightarrow \R
$$
by
$$
G(u)\equiv \chi_{R}(g(u))e^{-\frac{2}{3}f(u)}\,.
$$
We then have that $G(u)$ is the limit in measure, as $N\rightarrow\infty$, of
\begin{equation}\label{car}
\chi_{R}\Big(\|S_{N}(u)\|_{L^2(S^1)}^{2}-\alpha_{N}\Big)e^{-\frac{2}{3}\int_{S^1}(S_{N}u)(x)^{3}dx}.
\end{equation}
Indeed since $\chi_{R}(x)$ and $e^{-\frac{2}{3}x}$ are continuous real functions, we
have that $\chi_{R}(g_N(u))$ and $e^{-\frac{2}{3}f_{N}(u)}$ converge in the $\theta$ measure 
to $\chi_{R}(g(u))$ and $e^{-\frac{2}{3}f(u)}$ respectively. Then we use that the
convergence in measure is stable with respect to the product operation to
conclude that indeed (\ref{car}) converges to $G(u)$ in measure.
Thus
the function $G$ is measurable from $(H^{-\sigma}_{0}(S^1),{\mathcal B})$ to $\R$.
We are going to show that in fact $G\in L^{p}(H^{-\sigma}_{0}(S^1),{\mathcal B},d\theta)$ for all finite $p\geq 1$.
The main point is the following statement.
\begin{proposition}\label{main}
Let $1\leq p<\infty$. Then there exists $C>0$ such that for every $N\geq 1$,
$$
\Big\|\chi_{R}\big(\|S_{N}u\|_{L^2(S^1)}^{2}-\alpha_{N}\big)
e^{-\frac{2}{3}\int_{S^1}(S_{N}u)(x)^{3}dx}\Big\|_{L^p(d\theta(u))}\leq C\,.
$$
\end{proposition}
\begin{proof}
Our goal is to evaluate the function
$
\theta(A_{\lambda}),
$
where
$$
A_{\lambda}\equiv
\Big(
u\in H^{-\sigma}_{0}(S^1)\,:\,
\chi_{R}\big(\|S_{N}u\|_{L^2(S^1)}^{2}-\alpha_{N}\big)
e^{-\frac{2}{3}\int_{S^1}(S_{N}u)(x)^{3}dx}>\lambda
\Big)
$$
for $\lambda\geq 200$.
More precisely, we need to show the convergence and the uniform with respect
to $N$ boundedness of the integral $\int^{\infty}\lambda^{p-1}\theta(A_{\lambda})d\lambda$.
Set 
\begin{equation}\label{n0}
N_{0}\equiv (\log\lambda)^{2}\,.
\end{equation}
Suppose first that $N_0\geq N$. 
Using the H\"older inequality, we get for $u\in A_{\lambda}$,
\begin{eqnarray*}
\Big|\int_{S^1}(S_{N}u)(x)^{3}dx\Big|
& \leq & C\|S_N u\|_{L^2(S^1)}^{2}\|S_N u\|_{L^{\infty}(S^1)}
\leq  C\alpha_N \|S_N u\|_{L^{\infty}(S^1)}
\\
& \leq &  C\log(N) \|S_N u\|_{L^{\infty}(S^1)}
\leq C(\log\log\lambda)\|S_N u\|_{L^{\infty}(S^1)} \,.
\end{eqnarray*}
Hence for every $\delta>0$ there exist $C$ and $c$, independent of $N$, such that
$$
\theta(A_{\lambda})\leq 
\theta
\Big(
u\in H^{-\sigma}_{0}(S^1)\,:\, \|S_{N}u\|_{L^{\infty}(S^1)}\geq
c(\log\lambda)^{1-\delta},\,\,
\|S_{N}u\|_{L^2(S^1)}^{2}\leq C\log\log\lambda
\Big).
$$
Thus, using Proposition~\ref{calais} (with $(\log\lambda)^{1-\delta}$ instead
of $\lambda$), we infer that for every $\varepsilon>0$
there exist $C>0$, $c>0$ such that
$
\theta(A_{\lambda})\leq 
C\exp(-c(\log\lambda)^{2-\varepsilon})\leq C_{L}\lambda^{-L}
$
which yields the needed uniform integrability property.
\\

We can therefore suppose in the sequel of the proof that $N>N_0$, where $N_0$
is defined by (\ref{n0}).
Consider the set
$$
B_{\lambda,\kappa}\equiv\Big(u\in H^{-\sigma}_{0}(S^1)\,:\, |g_{N}(u)-g_{N_0}(u)|>\kappa\Big),
$$
where $g_{N}$ is defined by (\ref{gN}) and $\kappa$ is a large constant.
Lemma~\ref{Wickpak} yields
$$
\theta(B_{\lambda,\kappa})\leq Ce^{-\delta\kappa (\log\lambda)}= C\lambda^{-\delta\kappa}.
$$
Therefore if $\kappa\gg 1$ then $\mu(B_{\lambda,\kappa})\leq C\lambda^{-p-10}$.
Hence it suffices to evaluate $\theta(A_{\lambda}\backslash B_{\lambda,\kappa})$.
Let us observe that for $u\in A_{\lambda}\backslash B_{\lambda,\kappa}$ one has
\begin{eqnarray*}
\|S_{N_0}u\|_{L^2(S^1)}^{2} & =& (\|S_{N}u\|_{L^2(S^1)}^{2}-\alpha_{N})-
(g_{N}(u)-g_{N_0}(u))
+\alpha_{N_0}
\\
& \leq &  
C+\kappa+C\log(N_0)\leq C\log\log\lambda\,.
\end{eqnarray*}
Therefore $A_{\lambda}\backslash B_{\lambda,\kappa}\subset C_{\lambda}$ where
$$
C_{\lambda}\equiv
\Big(
u\in H^{-\sigma}_{0}(S^1)\,:\, 
\Big|
\int_{S^1}(S_{N}u)(x)^{3}dx
\Big|
\geq
\frac{3}{2}\log\lambda,\,\,
\|S_{N_0}u\|_{L^2(S^1)}^{2}\leq C\log\log\lambda
\Big).
$$
We next observe that
$C_{\lambda}\subset D_{\lambda}\cup E_{\lambda}$, where
$$
D_{\lambda}\equiv
\Big(
u\in H^{-\sigma}_{0}(S^1)\,:\, 
\Big|
\int_{S^1}(S_{N_0}u)(x)^{3}dx
\Big|
\geq
\frac{1}{2}\log\lambda,\,\,
\|S_{N_0}u\|_{L^2(S^1)}^{2}\leq C\log\log\lambda
\Big)
$$
and
$$
E_{\lambda}\equiv
\Big(
u\in H^{-\sigma}_{0}(S^1)\,:\, 
\Big|
\int_{S^1}(S_{N}u)(x)^{3}dx
-
\int_{S^1}(S_{N_0}u)(x)^{3}dx
\Big|
\geq
\frac{1}{2}\log\lambda
\Big).
$$
Using Lemma~\ref{distri} we obtain that
for every $\alpha<1/2$ there exists $C>0$ and $\delta>0$ such that
$
\theta(E_{\lambda})\leq
Ce^{-\delta (N_0^{\alpha}\log\lambda)^{2/3}}\leq C_{L}\lambda^{-L}
$
by taking $\alpha$ close enough to $1/2$
(recall that $N_0=(\log\lambda)^{2}$).
Hence it only remains to evaluate $\theta(D_{\lambda})$.
Using the H\"older inequality, we obtain that  for $u\in D_{\lambda}$ one has
$$
\Big|\int_{S^1}(S_{N_0}u)(x)^{3}dx\Big|
\leq 
\|S_{N_{0}}u\|_{L^{\infty}(S^1)}
\|S_{N_{0}}u\|_{L^2(S^1)}^{2}
\leq  
C\log\log\lambda\|S_{N_{0}}u\|_{L^{\infty}(S^1)}\,.
$$
Therefore for every $\delta>0$ there exists $C$ and $c$ such that
$$
\theta(D_{\lambda})\leq 
\theta
\Big(
u\in H^{-\sigma}_{0}(S^1)\,:\, \|S_{N_0}u\|_{L^{\infty}(S^1)}\geq
c(\log\lambda)^{1-\delta},\,\,
\|S_{N_0}u\|_{L^2(S^1)}^{2}\leq C\log\log\lambda
\Big).
$$
Using once again
Proposition~\ref{calais}, we infer that for every $\varepsilon>0$
there exist $C>0$, $c>0$ such that
$
\theta(D_{\lambda})\leq 
C\exp(-c(\log\lambda)^{2-\varepsilon})\leq C_{L}\lambda^{-L}.
$
Hence we conclude that
$$
\theta(A_{\lambda}\backslash B_{\lambda,\kappa})\leq\theta(C_{\lambda})\leq
\theta(D_{\lambda})+
\theta(E_{\lambda})\leq C_{L}\lambda^{-L}.
$$
This completes the proof of Proposition~\ref{main}.
\end{proof}
Let us now consider the sequence of measurable functions from
$(H^{-\sigma}_{0}(S^1),{\mathcal B})$ to $\R$ defined as
$$
G_{N}(u)\equiv 
\chi_{R}\big(\|S_{N}u\|_{L^2(S^1)}^{2}-\alpha_{N}\big)e^{-\frac{2}{3}\int_{S^1}(S_{N}u)(x)^{3}dx}\,.
$$
Since $G_{N}$ converges to $G$ in measure, we obtain that there exists a
subsequence $N_k$ such that 
$$
G(u)=\lim_{k\rightarrow \infty}G_{N_k}(u),\quad \theta\,\,\, {\rm a.s.}
$$
Proposition~\ref{main} implies that there exists a constant $C$ such that
$$
\|G_{N_k}(u)\|_{L^p(d\theta(u))}\leq C,\quad \forall k\in\N.
$$
Hence Fatou's lemma implies that $G(u)\in L^p(d\theta(u))$ and moreover 
$$
\int_{H^{-\sigma}_{0}(S^1)}|G(u)|^{p}d\theta(u)
\leq \liminf_{k\rightarrow\infty}
\int_{H^{-\sigma}_{0}(S^1)}|G_{N_k}(u)|^{p}d\theta(u)\,.
$$
Let now $h$ be a bounded continuous function from $H^{-\sigma}_{0}(S^1)$ to $\R$.
Our goal is to show that
\begin{equation}\label{krai}
\lim_{N\rightarrow \infty}\int_{H^{-\sigma}_{0}(S^1)}G_{N}(u)h(u)d\theta(u)=
\int_{H^{-\sigma}_{0}(S^1)}G(u)h(u)d\theta(u)\,.
\end{equation}
Let us fix $\varepsilon >0$. Consider the set
$$
A_{N,\varepsilon}\equiv \big(
u\in H^{-\sigma}_{0}(S^1)\,:\, |G_{N}(u)-G(u)|\leq \varepsilon
\big).
$$
Denote by $A_{N,\varepsilon}^{c}$ the complementary set in
$H^{-\sigma}_{0}(S^1)$ of $A_{N,\varepsilon}$. Then, using that $h$ is
bounded, Proposition~\ref{main}
and the Cauchy-Schwarz inequality, we infer that
$$
\Big|
\int_{A_{N,\varepsilon}^{c}}(G_{N}(u)-G(u))h(u)d\theta(u)
\Big|
\leq
C\|G_N-G\|_{L^2(d\theta)}[\theta(A_{N,\varepsilon}^{c})]^{\frac{1}{2}}\leq
C[\theta(A_{N,\varepsilon}^{c})]^{\frac{1}{2}}\,,
$$
where $C$ is independent of $N$ and $\varepsilon$.
On the other hand
$$
\Big|
\int_{A_{N,\varepsilon}}(G_{N}(u)-G(u))h(u)d\theta(u)
\Big|
\leq
C\varepsilon
$$
and thus we have (\ref{krai}) since the convergence in measure of $G_N$ to $G$
implies that for a fixed $\varepsilon$,
$$
\lim_{N\rightarrow \infty}\theta(A_{N,\varepsilon}^{c})=0\,.
$$
This completes the proof of Theorem~\ref{thm}.
\qed 
\\

Let us observe that for $R\gg 1$ the measure $d\rho(u)$ is not trivial. 
Indeed, by the estimates on second order Wiener chaos (see Lemma~\ref{Wickpak})
we infer that
$$
\theta\Big(u\in H^{-\sigma}_{0}(S^1)\,:\,\big|\|S_N u\|_{L^2(S^1)}^{2}-\alpha_{N}\big|>R\Big)\leq Ce^{-\delta R}\, ,
$$
where $C>0$ and $\delta>0$ are independent of $R$. 
Since for $R\geq 3$,
\begin{multline*}
\Big(u\in H^{-\sigma}_{0}(S^1)\,:\,|g(u)|>R\Big)
\subset
\Big(u\in H^{-\sigma}_{0}(S^1)\,:\,|g_{N}(u)|>R-2\Big)
\\
\cup
\Big(u\in H^{-\sigma}_{0}(S^1)\,:\,|g(u)-g_{N}(u)|>1\Big)
\end{multline*}
using the convergence in measure of $g_{N}(u)$ to $g(u)$, we obtain that
for $R\gg 1$ the set $(u\in H^{-\sigma}_{0}(S^1):|g(u)|\leq R)$ is of positive
$\theta$ measure and thus $d\rho(u)$ is a nontrivial measure since its (non-negative) density
is not vanishing on a set of positive $\theta$ measure.
The result of Theorem~\ref{thm} implies some additional properties of the
convergence of $\rho_N$ to $\rho$. For instance we have the following statement.
\begin{proposition}\label{liminf}
Let $U$ be an open set of $H^{-\sigma}(S^1)$. Then
\begin{equation}\label{open}
\liminf_{N\rightarrow\infty}\rho_{N}(U)\geq \rho(U).
\end{equation}
Let $V$ be a closed set of  $H^{-\sigma}(S^1)$. Then
\begin{equation}\label{closed}
\rho(V)\geq \limsup_{N\rightarrow\infty}\rho_{N}(V).
\end{equation}
\end{proposition}
\begin{proof}
Applying Theorem~\ref{thm} to $h=1$, we obtain that
\begin{equation}\label{gran}
\int_{H^{-\sigma}_{0}(S^1)}G(u)d\theta(u)
= \lim_{N\rightarrow\infty}
\int_{H^{-\sigma}_{0}(S^1)}G_{N}(u)d\theta(u)\,.
\end{equation}
We set
$$
\beta_{N}\equiv \int_{H^{-\sigma}_{0}(S^1)}G_{N}(u)d\theta(u),\quad
\beta\equiv \int_{H^{-\sigma}_{0}(S^1)}G(u)d\theta(u).
$$
If $\beta=0$ the assertion is trivial.
We can therefore suppose that $\beta\neq 0$ and that there exists $N_0$ such that $\beta_N\neq 0$, $\forall N\geq N_0$.
Next, we define the probability measures on $(H^{-\sigma}_{0}(S^1),{\mathcal B})$ as
$$
d\tilde{\rho}_{N}\equiv \beta_{N}^{-1}d\rho_{N},N\geq N_0,\quad d\tilde{\rho}\equiv \beta^{-1}d\rho.
$$
Since $\lim_{N\rightarrow\infty}\beta_{N}=\beta$ (see (\ref{gran}), Theorem~\ref{thm}
implies that
for every continuous bounded function $h$ from $H^{-\sigma}_{0}(S^1)$ to $\R$, we have
$$
\int_{H^{-\sigma}_{0}(S^1)}h(u)d\tilde{\rho}(u)=\lim_{N\rightarrow\infty}
\int_{H^{-\sigma}_{0}(S^1)}h(u)d\tilde{\rho}_{N}(u)\,.
$$
But it is known (see e.g. \cite{LQ,LT}) that the above convergence is in fact
equivalent with the fact that for every open set of $H^{-\sigma}(S^1)$ one has
\begin{equation}\label{openbis}
\liminf_{N\rightarrow\infty}\tilde{\rho}_{N}(U)\leq \tilde{\rho}(U).
\end{equation}
Using (\ref{gran}), we infer that (\ref{open}) holds.
Finally, one obtains (\ref{closed}) by passing to complementary sets in (\ref{openbis}).
This completes the proof of Proposition~\ref{liminf}.
\end{proof}
\section{The random distribution $\varphi(\omega,x)$ and the Benjamin-Ono equation}
\subsection{Behavior of the map $u\mapsto u^2$ on the statistical ensemble}
If one is interested to construct solutions of the Benjamin-Ono equation with
almost all data $\varphi(\omega,x)$  given by (\ref{stat}), in view of the
structure of the nonlinearity, it is natural to
ask about regularity properties of $\varphi^{2}(\omega,x)$. Since
$\|\varphi(\omega,\cdot)\|_{L^2}=\infty$ a.s. it is natural to project 
$\varphi^{2}(\omega,x)$ on the non zero modes. If we denote by $\Pi$ the
projector on the non zero modes, we have the following statement.
\begin{lemme}\label{carre}
For every $s<0$ there exists a constant $C$ such that for every $N$,
$$
\E\Big(\|\Pi(\varphi_{N}^{2}(\omega,x))\|^{2}_{H^{s}(S^1)}\big)\leq C. 
$$
\end{lemme}
\begin{remarque}
The nontrivial point is that $C$ is independent of $N$.
\end{remarque}
\begin{proof}
Write
$$
\Pi(\varphi_{N}^{2}(\omega,x))=
\sum_{\stackrel{0<|n_1|, |n_2|\leq N}{ n_1+n_2\neq 0}}
\frac{g_{n_1}(\omega)}{2\sqrt{\pi|n_1|}}\frac{g_{n_2}(\omega)}{2\sqrt{\pi|n_2|}}\, e^{i(n_1+n_2)x}\,.
$$
Therefore
$$
\Big\|
\Pi(\varphi_{N}^{2}(\omega,\cdot))
\Big\|_{H^s(S^1)}^{2}
=
\frac{1}{2}
\sum_{n\neq 0}\langle n\rangle^{2s}
\Big|
\sum_{\stackrel{
0<|n_1|, |n_2|\leq N
}{ n_1+n_2=n}}
\frac{g_{n_1}(\omega)}{\sqrt{|n_1|}}\frac{g_{n_2}(\omega)}{\sqrt{|n_2|}}
\Big|^{2}\,.
$$
Denote by $G(\omega)$ the right hand-side of the above equality. Then using
the independence of $g_{n}(\omega)$ one verifies that
$$
\E(G)\leq C
\sum_{n\neq 0}|n|^{2s}
\sum_{\stackrel{
0< |n_1|, |n_2|\leq N
}{ n_1+n_2=n}}
\frac{1}{|n_1|}
\frac{1}{|n_2|}
\leq
C\sum_{n\neq 0}\sum_{n_1\in\Z\backslash\{0,n\}}|n|^{2s}|n_1|^{-1}|n-n_1|^{-1}\,.
$$
Therefore, using Lemma~\ref{elem1}, we get
$$
\E(G)\leq C_{\varepsilon}\sum_{n\neq 0}|n|^{2s}(1+|n|)^{-1+\varepsilon}<\infty
$$
provided $2s+\varepsilon<0$, i.e. $0<\varepsilon<-2s$. This completes the proof
of Lemma~\ref{carre}.
\end{proof}
\subsection{Tao's gauge transform on the statistical ensemble}
In \cite{Tao}, Tao introduces a gauge transform which turns out to be a
crucial tool in the low regularity well-posedness of the Benjamin-Ono
equation (see \cite{Molinet,BP,IK}).
We now study the action of this transform on the functions on the
statistical ensemble (\ref{stat}). 
Recall that Tao's gauge transform is defined by
$$
u\longmapsto P_{+}\Big(e^{-i\partial_{x}^{-1}u}\, u\Big)\equiv \Phi(u),
$$
where $P_{+}$ denotes the projector on the positive frequencies.
For $u\in L^2(S^1)$ with zero mean value the gauge transform $\Phi(u)$ is
easily seen to belong to $L^2(S^1)$. For $u\in H^{s}_{0}(S^1)$, $-1/2<s<0$ one
may give sense of the product $e^{-i\partial_{x}^{-1}u}u$ in $L^1(S^1)$, since
$\partial_{x}^{-1}u\in H^{1-s}(S^1)$ a.s. and $1-s>1/2$ implies that 
$e^{-i\partial_{x}^{-1}u}\in H^{1-s}(S^1)$ a.s. It is however not a priori
clear that for $u\in H^{s}_{0}(S^1)$, $-1/2<s<0$, the transform $\Phi(u)$ is
also in $H^{s}_{0}(S^1)$. This turns out to be the case for
$u=\varphi(\omega,x)$ a.s. in $\omega$ as shows the next lemma.
\begin{lemme}\label{gauge}
Let us fix $s<0$. Then $\Phi(\varphi(\omega,\cdot))\in H^s_{0}(S^1)$ a.s.
\end{lemme}
\begin{proof}
Write
$$
e^{-i\partial_{x}^{-1}u}=\sum_{k=1}^{\infty}\frac{(-i\partial_{x}^{-1}u)^{k}}{k!}\,.
$$
Recall (see e.g. \cite{LQ}) that there exists an a.s. finite real valued
random variable $H(\omega)$ such that for every $n=1,2,\cdots$,
\begin{equation}\label{LQ}
|g_{n}(\omega)|\leq (\log(1+n))^{\frac{1}{2}}H(\omega).
\end{equation}
It suffices therefore to show that
\begin{equation}\label{rep}
\Big\|
P_{+}\Big((\partial_{x}^{-1}\varphi(\omega,\cdot))^{k}\varphi(\omega,\cdot)\Big)
\Big\|_{H^s(S^1)}
\leq
C\big(H(\omega)\big)^{k+1}\,.
\end{equation}
The proof of (\ref{rep}) is based on a repetitive use of Lemma~\ref{elem2}.
The square of the left hand-side of (\ref{rep}) can be bounded by
$$
C\sum_{n> 0}n^{2s}
\Big|
\sum_{\stackrel{n_{j}\neq 0,j=1,\cdots, k+1}{ n=n_1+\cdots+ n_{k+1}}}
\frac{g_{n_1}(\omega)}{|n_1|^{\frac{1}{2}}}
\frac{g_{n_2}(\omega)}{|n_2|^{\frac{3}{2}}}
\cdots
\frac{g_{n_{k+1}}(\omega)}{|n_{k+1}|^{\frac{3}{2}}}
\Big|^{2}\,.
$$
Next using $k+1$ times (\ref{LQ}) we bound the last expression as follows
$$
C_{\varepsilon}(H(\omega))^{2(k+1)}
\sum_{n> 0}n^{2s}
\Big|
\sum_{\stackrel{n_{j}\neq 0,j=1,\cdots, k+1}{ n=n_1+\cdots+ n_{k+1}}}
\frac{1}{|n_1|^{\frac{1}{2}-\varepsilon}}
\frac{1}{|n_2|^{\frac{3}{2}-\varepsilon}}
\cdots
\frac{1}{|n_{k+1}|^{\frac{3}{2}-\varepsilon}}
\Big|^{2}\,,
$$
where $\varepsilon\in]0,\frac{1}{4}[$.
Using Lemma~\ref{elem2}, we bound the above expression by
\begin{multline*}
C_{\varepsilon}(H(\omega))^{2(k+1)}
\\\times
\sum_{n> 0}n^{2s}
\Big|
\sum_{\stackrel{n_{j}\neq 0,j=2,\cdots, k}{ n\neq n_2+\cdots+ n_{k}}}
\frac{1}{|n_2|^{\frac{3}{2}-\varepsilon}}
\cdots
\frac{1}{|n_{k}|^{\frac{3}{2}-\varepsilon}}
\frac{1}{|n-n_2-\cdots-n_{k}|^{\frac{1}{2}-\varepsilon}}
\Big|^{2}.
\end{multline*}
Finally using $k-1$ more times Lemma~\ref{elem2},
we eliminate consequently $n_k$, $n_{k-1}$ etc. up to $n_2$
and thus we bound the last expression
by
$$
C_{\varepsilon}(H(\omega))^{2(k+1)}
\sum_{n> 0}n^{2s}
\Big|
\frac{1}{|n|^{\frac{1}{2}-\varepsilon}}
\Big|^{2}\leq C_{\varepsilon}(H(\omega))^{2(k+1)},
$$
provided $0<\varepsilon<\min(\frac{1}{4},-s)$. This proves (\ref{rep}) and Lemma~\ref{gauge}
is therefore established.
\end{proof}
\subsection{Bounds on the second Picard  iteration associated to the Benjamin-Ono equation with data  $\varphi(\omega,x)$}
If we set $\sigma(n)\equiv- n|n|$, we then have that for $w\in H^s_0(S^1)$,
$s\in\R$ the solution of the linearized around the zero solution Benjamin-Ono
equation
$$
(\partial_t+H\partial_{x}^2)u=0,\quad u|_{t=0}=w
$$
is given by
$$
u(t,x)\equiv \exp(-tH\partial_x^2)(w)=\sum_{n\neq 0}
e^{it\sigma(n)}e^{inx}\hat{w}(n)\,.
$$
If we are interested to solve the  Benjamin-Ono
equation with initial data $\varphi(\omega,x)$ 
then it is useful to study the problem
\begin{equation}\label{second} 
(\partial_t+H\partial_{x}^2)u+\partial_{x}
\Big(\exp(-tH\partial_x^2)(\varphi(\omega,\cdot))\Big)^{2},\quad u|_{t=0}=0\,.
\end{equation}
If we replace in (\ref{second}) $\varphi(\omega,x)$ by an $H^s_0(S^1)$ ,
$s\geq 0$ function then
it follows from the work by Molinet \cite{Molinet} that the solution of (\ref{second}) is in
$H^s_0(S^1)$. 
We are now going to show that in the context of (\ref{second}) the solution is
a.s. in all $H^s_{0}(S^1)$, $s<0$.
\begin{proposition}\label{kraino}
For every $s<0$, the solution $u$ of (\ref{second}) is a.s. in $H^s(S^1)$.
\end{proposition}
\begin{proof}
We have that
\begin{multline*}
\partial_{x}\Big(\exp(-tH\partial_x^2)(\varphi(\omega,\cdot))\Big)^{2}
\\
=
\sum_{n_1\neq 0,n_2\neq 0}
i(n_1+n_2)
\frac{g_{n_1}(\omega)}{2\sqrt{\pi |n_1|}}
\frac{g_{n_2}(\omega)}{2\sqrt{\pi |n_2|}}
e^{it(\sigma(n_1)+\sigma(n_2))}e^{i(n_1+n_2)x}\,.
\end{multline*}
On the other hand by the Duhamel principle the solution of (\ref{second}) is
given by
$$
u(t,x,\omega)
=
-\int_{0}^{t}
\exp(-(t-\tau)H\partial_x^2)
\Big(
\partial_{x}\Big(\exp(-\tau H\partial_x^2)(\varphi(\omega,\cdot))\Big)^{2}
\Big)
d\tau.
$$
Therefore there exists a numerical constant $c$ such that
$$
u(t,x,\omega)
=
c\sum_{n\neq 0}n e^{it\sigma(n)}\Big(
\sum_{n_1\neq 0,n}
\frac{e^{it(\sigma(n_1)+\sigma(n-n_1)-\sigma(n))}-1}{\sigma(n_1)+\sigma(n-n_1)-\sigma(n)}
\,\,
\frac{g_{n_1}(\omega)}{\sqrt{|n_1|}}
\frac{g_{n-n_1}(\omega)}{\sqrt{|n-n_1|}}
\Big)e^{inx}\,.
$$
Using a direct case by case analysis implies that for $n_1\neq 0,n$, $n\neq0$,
$$
|\sigma(n_1)+\sigma(n-n_1)-\sigma(n)|\geq |n|.
$$
Therefore using the independence of $g_n(\omega)$ and Lemma~\ref{elem1}, we
obtain for $s<0$
$$
\|u(t,\cdot,\cdot)\|_{L^{2}(\Omega;H^s_{0}(S^1))}^{2}
\leq
C
\sum_{n\neq 0}\sum_{n_1\neq 0,n}\frac{|n|^{2s}}{|n_1 (n-n_1)|}
\leq
\sum_{n\neq 0}
\frac{C_{\varepsilon}|n|^{2s}}{|n|^{1-\varepsilon}}<\infty,
$$
provided $\varepsilon>0$ being such that $2s+\varepsilon<0$.
This completes the proof of Proposition~\ref{kraino}.
\end{proof}

\end{document}